\documentclass[a4paper,fleqn]{cas-dc}

\usepackage[numbers,sort&compress]{natbib}
\usepackage{enumitem}
\usepackage{varwidth}
\usepackage{subcaption}
\usepackage{amsmath}
\usepackage{mathtools}
\def\tsc#1{\csdef{#1}{\textsc{\lowercase{#1}}\xspace}}
\tsc{WGM}
\tsc{QE}

\begin{document}
\setcitestyle{square,numbers,sort&compress} 
\bibpunct{[}{]}{,}{n}{}{,}
\setlength{\mathindent}{0pt}
\let\WriteBookmarks\relax
\def\floatpagepagefraction{1}
\def\textpagefraction{.001}
\shorttitle{Parity in Cyclic Communities}    

\shortauthors{T.M. Murdie et al.}  

\title [mode = title]{Effect of Parity in Cyclically Competing Communities}  



%

\author[]{T.M. Murdie}[orcid=0009-0001-0112-9848]

\cormark[1]

\ead{tenacitymurdie2030@u.northwestern.edu}


\credit{Writing - review \& editing, Writing - original draft, Methodology, Investigation, Formal analysis, Conceptualization}

\affiliation[]{organization={Department of Engineering Sciences and Applied Mathematics},
            addressline={Northwestern University}, 
            city={Evanston},
            postcode={60208}, 
            state={IL},
            country={USA}}

\author[]{D.M. Abrams}[orcid=0000-0002-6015-8358]

\ead{dmabrams@northwestern.edu}


\credit{Writing - review \& editing, Writing - original draft, Methodology, Investigation, Formal analysis, Conceptualization}

\author[]{A. Bayliss}

\ead{a-bayliss@northwestern.edu}
\credit{Writing - review \& editing, Writing - original draft, Methodology, Investigation, Formal analysis, Conceptualization}

\author[]{V.A. Volpert}
\ead{v-volpert@northwestern.edu}
\credit{Writing - review \& editing, Writing - original draft, Methodology, Investigation, Formal analysis, Conceptualization}

\cortext[1]{Corresponding author}



\begin{abstract} 
Communities exhibiting cyclic (i.e., ``winnerless'') competition occur throughout nature. These systems can exhibit qualitatively different long-term behavior depending on the strength of inter-species competition and a specific property of the number of interacting species: parity. We consider how ecological communities can adapt to and transition between an odd and an even number of cyclically competing species. Previous studies have shown that, for strong inter-species competition, odd parity communities are dynamically unstable, while species in even parity communities form stable alliances of maximal noncompeting sets. We trace a homotopy between two May-Leonard type models: the odd parity $N=3$
model and the even parity $N=4$ system. \\
\indent{}We identify all physical steady states and analyze their stabilities. We show the existence of stable transitional states that are unique to our symmetric model. We then use a piecewise linear approximation technique to characterize the stability of heteroclinic cycles. We use numerical simulations to map distinct regimes in parameter space. Our computations confirm our analytical results and also reveal a window in parameter space where limit cycles occur. We illustrate that parity transitions can be complex and characterized by diverse parameter-dependent pathways.
\end{abstract}

\begin{graphicalabstract}
\includegraphics[width=13cm, height=5cm]{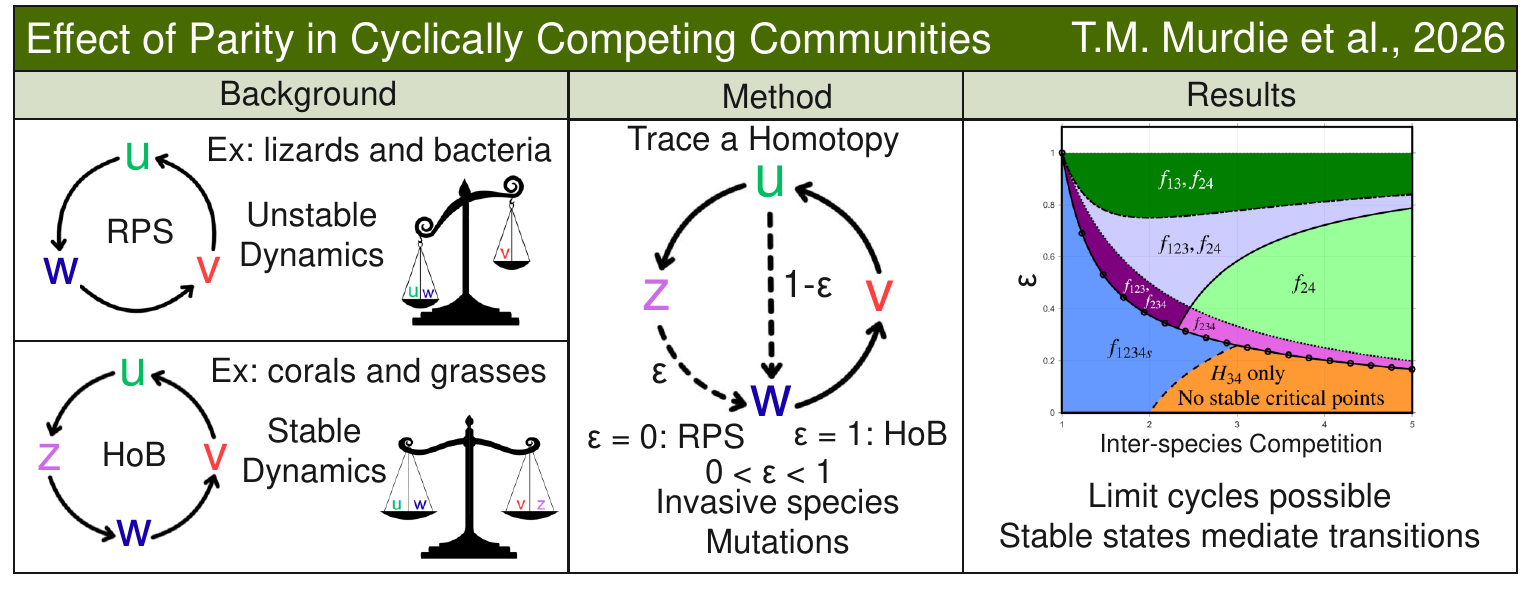}
\end{graphicalabstract}

\begin{highlights}
\item Parity leads to qualitatively different behavior in systems with cyclic competition.
\item Altered species interactions can break alliances and cause unstable system dynamics. 
\item Unique stable states mediate the transition from three to four competing species.
\item Numerical simulations identified a window of stable limit cycles in parameter space.
\end{highlights}


\begin{keywords}
 Population dynamics\sep Cyclic ecosystems\sep Species interactions\sep Heteroclinic cycles\sep Limit cycles
\end{keywords}

\maketitle

\section{Introduction}



\begin{figure}
  \centering
    \includegraphics[width=1\linewidth, height=8cm, keepaspectratio, trim={10 25 10 10}, clip]{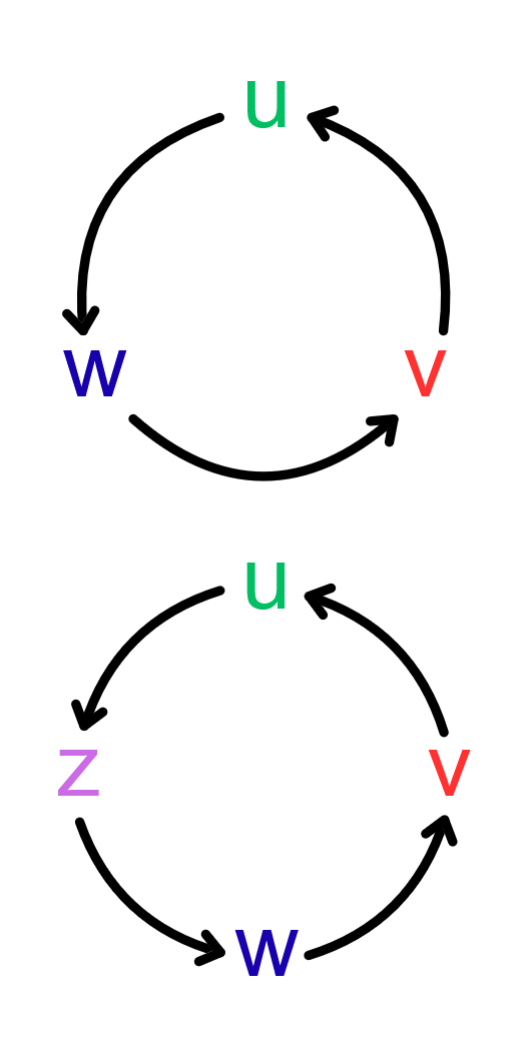}
    \caption{Competition schemes for $N=3$ (top) and $N=4$ (bottom). Arrows point toward the species at a competitive disadvantage.}\label{fig:RPS_HoB}
\end{figure}

\indent{}\indent{}A simple solution to settle an argument between children is to have them play rock-paper-scissors, a game in which rock destroys scissors, scissors cut paper, and paper crumples around rock. No single strategy is guaranteed to win, ensuring a fair resolution to childish disputes. This type of intransitive competition is also present in the natural world. Non-hierarchical cycles of interacting species can promote biodiversity in bacteria \cite{kerr2002local,kirkup2004antibiotic,czaran2002chemical}, explain the oscillating prevalence of multiple male reproductive strategies in side-blotched lizards \cite{sinervo1996rock}, and shape the community structure of coral reefs and grasslands \cite{buss1979competitive, silvertown1994spatial, cameron2009parasite}.\\
\indent{}Although these examples portray ``winnerless'' competition such that each species wins against one neighbor in the cycle and loses to another, the absence of a dominating strategy does not inherently promise coexistence. Certainly, sufficiently weak inter-species competition results in the stable coexistence of all species in any May-Leonard-type model with $N$ species \cite{may1975nonlinear}, but more interesting results arise when this competition is strong. In fact, strong inter-species competition can enable one or more species to monopolize resources and drive competitors to near-extinction. 
In this paper, we focus our attention on communities with $N=3$ and $N=4$ cyclically competing species (Fig. \ref{fig:RPS_HoB}).

\begin{figure}[]
    \centering
    \begin{subfigure}[t]{0.45\textwidth}
        \centering
        \includegraphics[width=1\linewidth]{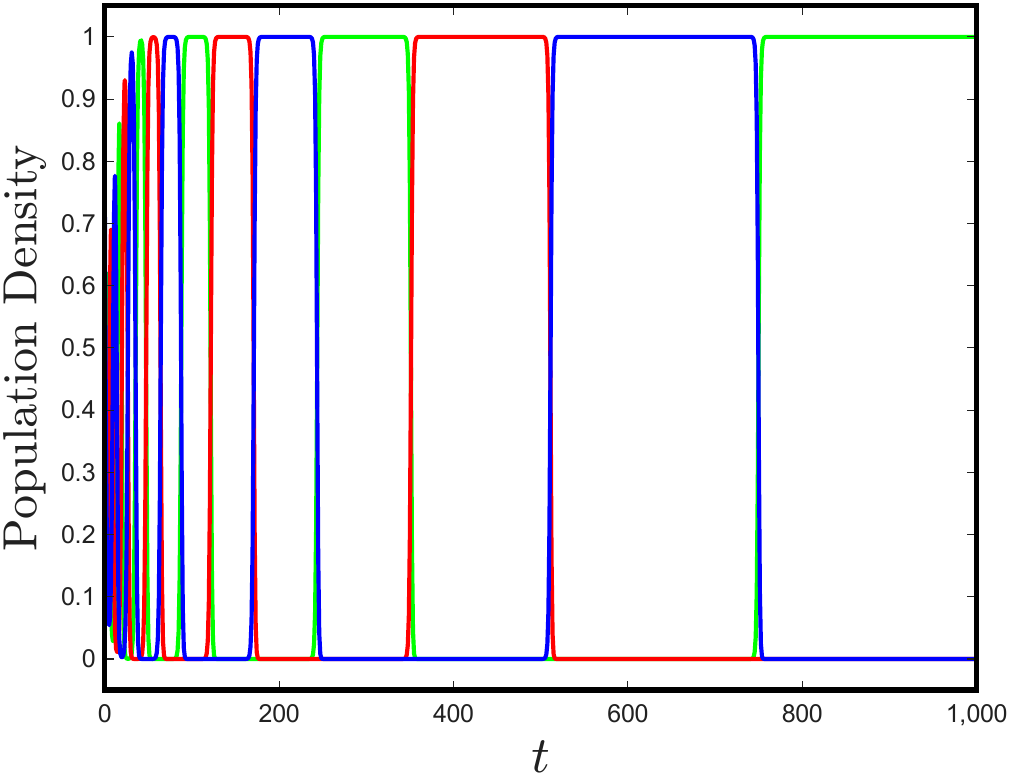}
        \caption{Heteroclinic cycle solution in a community with $N=3$ species, showing population density for species $u$ (green), $v$ (red), and $w$ (blue).\\}
        \label{subfig:RPS_Hcycle}
    \end{subfigure}
    \vfill
    \begin{subfigure}[t]{0.45\textwidth}
        \centering
        \includegraphics[width=0.9\linewidth]{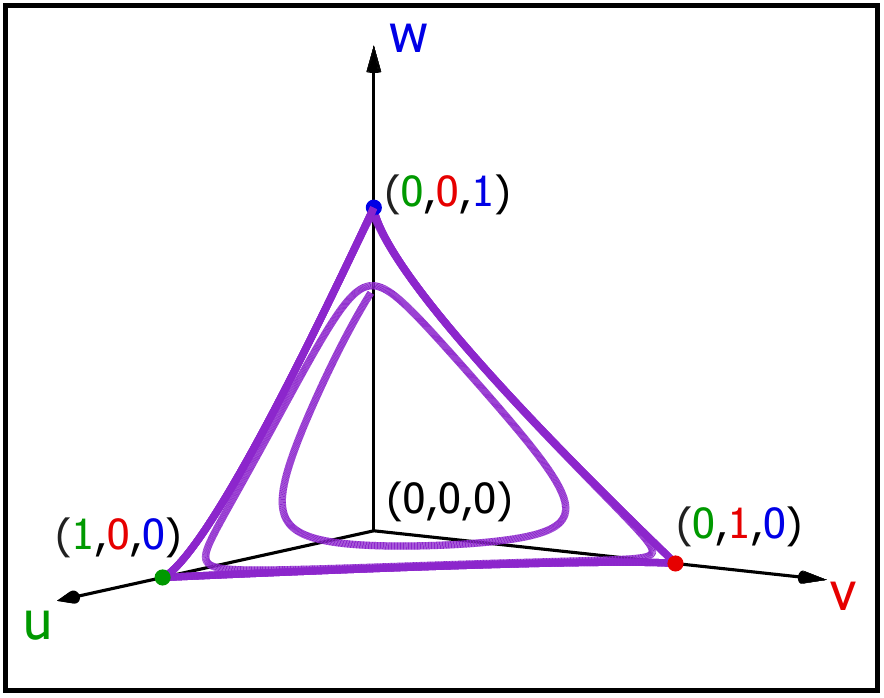}
        \caption{Three-dimensional population space, where the trajectory is shown visiting each single-species saddle.}
        \label{subfig:RPS_phase_space}
    \end{subfigure}
    \vfill
    \caption{Dynamic instability characteristic of odd parity communities with strong inter-species cyclic competition.}
    \label{fig:RPS_example}
\end{figure}

When $N=3$, the community is dynamically unstable: each species successively out-competes the others before its reign is replaced by that of its cyclic successor (Fig. \ref{subfig:RPS_Hcycle}). Mathematically, these alternating periods of dominance occur because the single-species fixed points are saddles with two stable and one unstable direction. The solution trajectory visits each saddle point in phase space, increasingly slowing down and spending more time near each saddle with every revolution (Fig. \ref{subfig:RPS_phase_space}). The closed connection of saddle points by a solution trajectory forms a heteroclinic cycle \cite{guckenheimer1988structurally}, a phenomenon observed---but not formally named---in ecological models by May and Leonard \cite{may1975nonlinear}. This dynamic instability persists in communities of odd $N>3$, though the heteroclinic cycle now connects saddles of partial alliances of $\frac{(N-1)}{2}$ mutually noncompeting species \cite{bayliss2020beyond}.\\

\indent{}The long-term behavior of a system is qualitatively different when $N$ is even. Though alliances of non-competing species still form, the system is now dynamically stable. That is, each alliance is now a stable equilibrium. The system exhibits bi-stability rather than cycles, and the winning alliance is determined by the initial conditions. The $N=4$ case has been described as the ``Hand of Bridge'' model, because the bi-stability between two pairs of species is reminiscent of pairs in a game of bridge \cite{durney2012stochastic}.\\
\indent{}There exists extensive literature pertaining to models with cyclic competition \cite{szolnoki2014cyclic}, and the effects of parity in stochastic systems with an arbitrary number of cyclically competing species have been characterized \cite{durney2011saddles}. However, a comprehensive understanding of how a system adapts to changes in parity is significantly limited. Mutations, invasive species, or environmental changes may alter the strength and number of species interactions \cite{bohannan2000linking,tuomainen2011behavioural}. It is our objective to understand how community behavior may be altered by changes in parity. In this paper, we modify the basic model introduced by May and Leonard \cite{may1975nonlinear} to allow for a homotopy between an odd parity rock-paper-scissors (RPS) community and an even parity Hand of Bridge community. We analyze the system in terms of stable and physical equilibrium states, as well as the presence and stability of heteroclinic cycles. Lastly, we present numerical simulations for different regimes observed across the transition.

\section{Model}

\begin{figure}[]
  \centering
    \includegraphics[width=1\linewidth, height=4cm, keepaspectratio, trim={10 255 10 10}, clip]{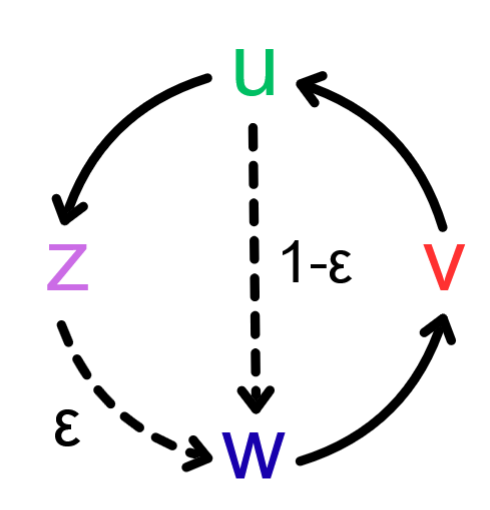}
    \caption{Competition scheme for the $3 \rightarrow 4$ species transition. Arrows point toward the species at a competitive disadvantage. Dashed lines indicate the competition $w$ experiences from $z$ (weighted by $\varepsilon$) and $u$ (weighted by $1-\varepsilon$).}\label{fig:3to4_model}
\end{figure}

Our model considers the transition between an RPS community and a Hand of Bridge community and traces homotopies between the two via a weighting parameter $0\le\varepsilon \le 1$ (Fig. \ref{fig:3to4_model}). 

We assume a well mixed population and use the symbols $u,v,w,$ and $z$ to represent the population densities of each species participating in cyclic competition. A suitable nondimensionalization results in the following system of equations:
\begin{subequations} \label{eqs:3to4_system}
    \begin{align}
    &\frac{du}{dt} = u(1-u-\alpha v),\\
    &\frac{dv}{dt} = v(1-v-\alpha w),\\
    &\frac{dw}{dt} = w(1-w-\varepsilon\alpha z - (1-\varepsilon)\alpha u),\\
    &\frac{dz}{dt} = z (1-z-\alpha u).
\end{align}
\end{subequations}
Here, $\alpha$ denotes the inter-species competition coefficient. As in \cite{may1975nonlinear}, we consider intra-species crowding but neglect crowding between different species, instead focusing on cyclic inter-species competition. For simplicity, we have assumed equal birth rates and competition coefficients for each species. Each species has a carrying capacity nondimensionalized to unity. \\
\indent{}When $\varepsilon=0$, our model describes an RPS community with an additional, but noncompetitive, species ($z$) present. That is, $w$ has an advantage over $v$, $v$ has an advantage over $u$, $u$ has an advantage over $w$ (thus forming the RPS cycle), $u$ also has an advantage over $z$, but for $\varepsilon=0$, $z$ does not exert competitive pressure on any other species in the community. As in standard RPS communities, the coexistence critical point becomes unstable for $\alpha>2$ due to an oscillatory instability, and an attracting heteroclinic cycle connects the single-species saddle points.\\
\indent{}In contrast, when $\varepsilon=1$, our model describes a Hand of Bridge community. In this system, alliances of non-directly competing species (i.e., the $u\text{-}w$ alliance and the $v\text{-}z$ alliance) are stable for $\alpha>1$. We note that both the RPS system and the Hand of Bridge system experience stable coexistence when $\alpha<1$. Thus, we limit the study of our $3 \rightarrow 4$ species model to $\alpha \ge 1$, where the effects of parity begin to emerge.

\section{Critical Points}

We label each critical point of (\ref{eqs:3to4_system}) with the prefix $f$ and with the indices of the non-zero components of the solution vector $(u,v,w,z)$ as subscripts. For example, $f_{1234}$ is the fixed point where all species coexist, while $f_{24}$ represents the steady state of the $v\text{-}z$ alliance.\\

\noindent{}\textbf{All components zero.} The extinction state $f_{0} = (0,0,0,0)$ is always unstable.\\

\noindent{}\textbf{Single-species states.} The states 
\begin{align*}
    &f_1 = (1,0,0,0),\ f_2 = (0,1,0,0),\\
    &f_3= (0,0,1,0), \text{ and } f_4=(0,0,0,1)
\end{align*}
are all saddle points with at least one positive and one negative eigenvalue of the associated Jacobian.\\

\noindent{}\textbf{Two-species states.} The states 
\begin{align*}
    &f_{12}=(1-\alpha,1,0,0),\ f_{14}=(1,0,0,1-\alpha),\\ &\text{and } f_{23}=(0,1-\alpha,1,0)\ 
\end{align*} 
are not physical for $\alpha>1$ and will not be considered further.\\

\noindent{}The state 
\begin{align*}
    f_{34}=(0,0,1-\varepsilon\alpha,1)
\end{align*}
is physical only for 
\begin{align*}
    \varepsilon<\varepsilon_{34} \text{, where } \varepsilon_{34} \coloneqq \frac{1}{\alpha},
\end{align*} 
and is always a saddle in this regime. The $v\text{-}z$ alliance state 
\begin{align*}
    f_{24}=(0,1,0,1)
\end{align*} 
is clearly always physical and is stable when $\alpha>1$ and 
\begin{align*}
    \varepsilon>\frac{1}{\alpha}=\varepsilon_{34}.
\end{align*} 
Lastly, the $u\text{-}w$ alliance state
\begin{align*}
    f_{13} = (1,0,1-(1-\varepsilon)\alpha,0)
\end{align*}
is physical when 
\begin{align*}
    \varepsilon>1-\frac{1}{\alpha}
\end{align*} 
and is stable when 
\begin{align*}
    \varepsilon>\varepsilon_{13} \text{, where } \varepsilon_{13}\coloneqq 1-\frac{1}{\alpha}+\frac{1}{\alpha^2}.
\end{align*}

\noindent{}\textbf{Three-species states.} The states
\begin{align*}
    &f_{124}=(1-\alpha,1,0,1-\alpha+\alpha^2) \\ &\text{and } f_{134} = (1,0,1-\alpha+\varepsilon\alpha^2,1-\alpha)\ 
\end{align*} 
are not physical when $\alpha>1$.\\

\noindent{}The state
\begin{align*}
    f_{234} = (0, 1-\alpha+\varepsilon\alpha^2, 1-\varepsilon\alpha,1)
\end{align*} 
is physical for 
\begin{align*}
    \frac{1}{\alpha}-\frac{1}{\alpha^2}<\varepsilon<\frac{1}{\alpha}=\varepsilon_{34}.
\end{align*}
We find that $f_{234}$ is stable for 
\begin{align*}
    \varepsilon_{**}<\varepsilon<\varepsilon_{34},
\end{align*}
where 
\begin{align*}
    \varepsilon_{**}\coloneqq \frac{1}{\alpha}-\frac{1}{\alpha^2}+\frac{1}{\alpha^3}.
\end{align*}
We refer the reader to Appendix \ref{apdx:f234} for the Jacobian matrix and corresponding eigenvalues.\\

\noindent{}The state 
\begin{align*}
    f_{123} = \Big( & \frac{1-\alpha+\alpha^2}{\alpha^3(1-\varepsilon)+1}, 
    \frac{1-\alpha+\alpha^2(1-\varepsilon)}{\alpha^3(1-\varepsilon)+1}, \\
    & \frac{1+(1-\varepsilon)(\alpha^2-\alpha)}{\alpha^3(1-\varepsilon)+1}, 0 \Big)
\end{align*}
is physical when 
\begin{align*}
    \varepsilon<1-\frac{1}{\alpha}+\frac{1}{\alpha^2}=\varepsilon_{13}.
\end{align*}
It is shown in Appendix \ref{apdx:f123} that $f_{123}$ is stable when the following inequality is satisfied:
\begin{align*}
    \text{max}(\varepsilon_{123},\varepsilon_{**}) < \varepsilon<\varepsilon_{13},
\end{align*}
where 
\begin{align*}
    \varepsilon_{123}\coloneqq \frac{(\alpha-1)^3-1}{\alpha(\alpha-1)^2}.
\end{align*}
For $\varepsilon=0$, $f_{123}$ reduces to the three-species coexistence state for the RPS model:
\begin{align*}
    f_{123} = \left(\frac{1}{\alpha+1},\frac{1}{\alpha+1},\frac{1}{\alpha+1},0 \right),
\end{align*}
such that species $z$ is completely extinct and only the RPS loop remains. In the three-species RPS system, this coexistence state is stable for $\alpha<2$. However, in our modified system, when $\varepsilon=0$, $f_{123}$ is a saddle for all $\alpha$. Although $\varepsilon=0$ results in species $z$ decoupling from the network of cyclic competition, it is impossible for species $z$ to go extinct while species $w$ survives if both $z$ and $w$ grow at the same rate and experience the same competitive pressure from $u$. That is, even though $f_{123}$ reduces to the RPS coexistence state for $\varepsilon=0$, our model features the stable four-species coexistence state $f_{1234s}$ (discussed below).\footnotemark[1]\\

\footnotetext[1]{Thus, for $\varepsilon=0$, the state $f_{123}$ becomes physically, though not mathematically, equivalent to the three-species coexistence state in the RPS system with symmetric competition.}

\noindent{}\textbf{Four-species Coexistence.} There are two different four-species coexistence steady states in this system: $f_{1234s}$ and $f_{1234a}$.\\

\noindent{}The coexistence state 
\begin{align*}
    f_{1234s}=\left(\frac{1}{\alpha+1},\frac{1}{\alpha+1},\frac{1}{\alpha+1},\frac{1}{\alpha+1} \right)
\end{align*}
is symmetric ($u^*=v^*=w^*=z^*=\frac{1}{\alpha+1}$). As shown in Appendix \ref{apdx:f1234s}, this symmetric coexistence state is stable in the regime
\begin{align*}
    \varepsilon_*=1-\frac{4}{\alpha}+\frac{8}{\alpha^2}-\frac{8}{\alpha^3}<\varepsilon<\frac{1}{\alpha}-\frac{1}{\alpha^2}+\frac{1}{\alpha^3}=\varepsilon_{**}.
\end{align*} 
This condition is only realizable when $\alpha \le 3$.\\

At the point $\varepsilon=\varepsilon_{**}$, there is a line of asymmetric coexistence states 
\begin{align*}
    f_{1234a}=\left(U,\frac{1-U}{\alpha},\frac{\alpha-1+U}{\alpha^2},1-\alpha U \right),
\end{align*}
where $U$ is a free parameter. More information is provided in Appendix \ref{apdx:f1234a}.\\

\begin{figure}[]
    \centering
    \includegraphics[width=1\linewidth]{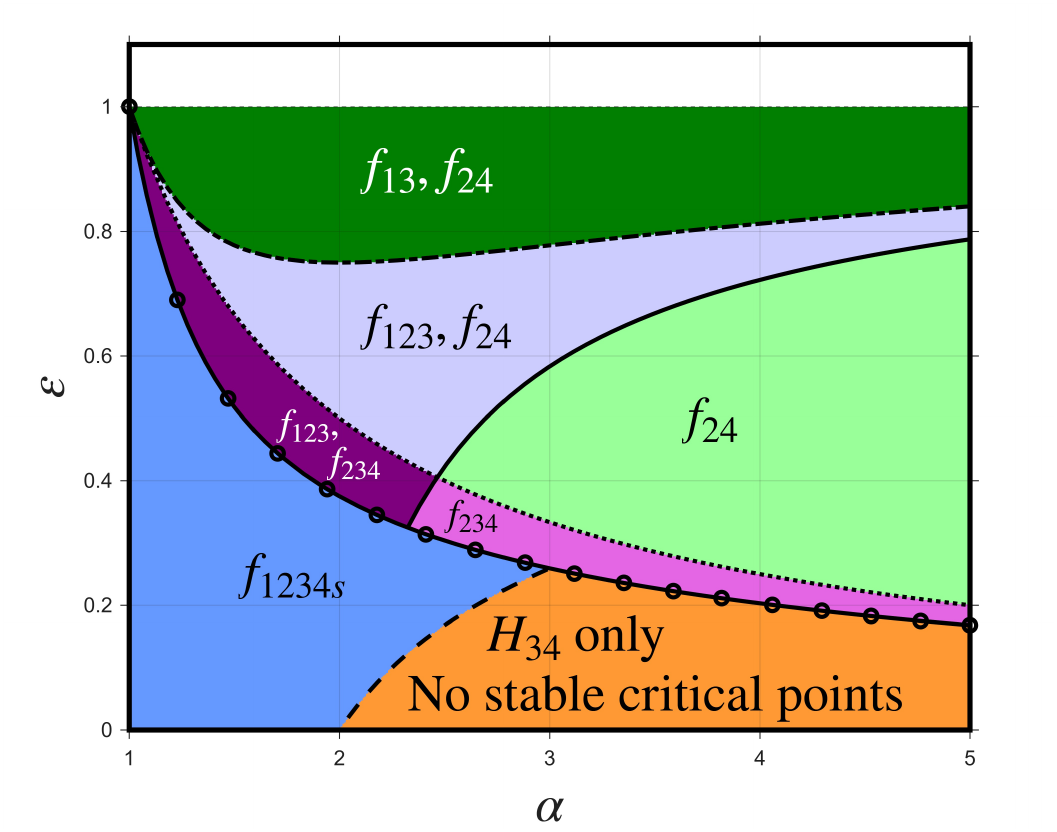}
    \caption{Regions of steady state stability in the $\alpha\text{-}\varepsilon$ plane, where $\alpha$ represents inter-species competition and $\varepsilon$ is the homotopy parameter. The colored regions indicate the following states are stable: $f_{1234s}$ (blue), $f_{123}$ (dark purple and lilac), $f_{234}$ (dark purple and magenta), $f_{13}$ (dark green), and $f_{24}$ (light green, lilac, and dark green). All critical points are unstable in the orange region. The critical $\varepsilon$ values are as follows: $\varepsilon_{**}$ (solid line with circles, where the system obtains the curve of steady states $f_{1234a}$), $\varepsilon_{*}$ (dashed line), $\varepsilon_{34}$ (dotted line), $\varepsilon_{123}$ (solid line), and $\varepsilon_{13}$ (dashed-dotted line). See figure \ref{fig:Regions_of_HCycle_Stability_3to4Model} for details regarding the stability of heteroclinic cycles.}
    \label{fig:Regions_of_Stability_3to4Model}
\end{figure}

\noindent{}\textbf{Summary.} Through stability analysis of the steady states, we can begin to trace the homotopy between the three-species RPS model and the four-species Hand of Bridge model. 
We find that the three-species state $f_{123}$ is physical for $\varepsilon=0$ but not for $\varepsilon=1$, while $f_{234}$ is non-physical at both boundary values of $\varepsilon$. Furthermore, at the extremes where $\varepsilon$ equals zero or unity, both three-species states are unstable. Instead, the states $f_{123}$ and $f_{234}$ are stable for intermediate values of $\varepsilon$ and mediate the transition from three to four cyclically competing species (odd to even parity). 
For fixed competition coefficient within the range $1<\alpha\le2$, stability transitions sequentially from the symmetric coexistence steady state $f_{1234s}$ to three-species and then two-species coexistence states as $\varepsilon$ increases (Fig. \ref{fig:Regions_of_Stability_3to4Model}). For fixed $\alpha>2$, the system transitions from unstable RPS-like dynamics with cycles of competitive exclusion (to be discussed later) to transitional coexistence steady states to the dynamically stable Hand of Bridge alliances of noncompeting pairs as $\varepsilon$ increases from zero to unity (Fig. \ref{fig:Regions_of_Stability_3to4Model}).

\section{Heteroclinic Cycles}

The hyperplanes $u=0$, $v=0$, $w=0$, and $z=0$ are each invariant manifolds for the system (\ref{eqs:3to4_system}). The heteroclinic orbits connecting two saddle points (e.g., $f_1 \rightarrow f_2$) lie in the corresponding hyperplanes and a connected cycle of these separatrices forms a heteroclinic cycle. Through numerical simulation, we identified the existence of three attracting heteroclinic cycles: $f_1 \rightarrow f_2 \rightarrow f_3 \rightarrow f_{34} \rightarrow f_1$ ($H_{34}$), $f_1 \rightarrow f_2 \rightarrow f_3 \rightarrow f_{4} \rightarrow f_1$ ($H_4$), and $f_1 \rightarrow f_2 \rightarrow f_3  \rightarrow f_1$ ($H_3$).
To examine the stability of a heteroclinic cycle, consisting of saddle points linked by separatrices, we pick an arbitrary point close to the cycle and compute the trajectory connecting the saddles by solving the piecewise linear systems obtained by linearization about the relevant saddle points. With the notable exception of $H_4$, we conclude that a cycle is stable if the trajectory gets closer to the cycle after one revolution from the initial point. This Poincar\'e mapping technique has been successfully implemented in previous papers \cite{bayliss2020beyond,bayliss2019mathematical} and is an application of more rigorous methods employed by \cite{krupa1995asymptotic,krupa2004asymptotic,melbourne1989heteroclinic,postlethwaite2010resonance}.\\

\noindent{}$\mathbf{H_{34}\textbf{: }f_1 \rightarrow f_2 \rightarrow f_3 \rightarrow f_{34} \rightarrow f_1.}$ \\\noindent{}\textbf{Step 1:} The initial point 
\begin{align*}
    B_1=(u_1, v_1, w_1,z_1)
\end{align*}
lies close to the separatrix $f_{34} \rightarrow f_1$, which exists in the hyperplane where $v=0$. Thus, $v_1 = \delta>0$ is small and $u_1,w_1,$ and $z_1$ are $\mathcal{O}(1)$ quantities. Consequently, we define the point $B_1$ to be
\begin{align*}
        B_1=(\mathcal{O}(1), \delta, \mathcal{O}(1),\mathcal{O}(1)).
\end{align*}
When the trajectory travels past $f_1$, the subsequent point $B_2$ will be between the saddles $f_1$ and $f_2$, where the $w$ and $z$ coordinates of both saddles are zero. Therefore, we expect 
\begin{align*}
    B_2=(\mathcal{O}(1),\mathcal{O}(1),\text{small}, \text{small}).
\end{align*}
We note that our analysis is only meaningful when $\delta\ll 1$, so that the $\mathcal{O}(1)$ terms may be safely neglected. By linearizing (\ref{eqs:3to4_system}) about $f_1$ and solving the resulting system, we find 
\begin{align*}
    &B_1=(\mathcal{O}(1), \delta, \mathcal{O}(1),\mathcal{O}(1)) \rightarrow \\
    &B_2 = (\mathcal{O}(1),\mathcal{O}(1),\mathcal{O}(1) \delta^{(1-\varepsilon)\alpha-1},\mathcal{O}(1)\delta^{\alpha-1}).
\end{align*}
For $w_2$ to be small, such that the trajectory remains close to the heteroclinic cycle, we require 
\begin{align}\label{eqs:H4_step1_condition}
    (1-\varepsilon)\alpha>1\implies \varepsilon<1-\frac{1}{\alpha}.
\end{align}
Likewise, $\alpha>1$ ensures $z_2$ is small. The full details of the linearization are provided in Appendix \ref{apdx:H34}. Subsequent steps rely on the same methodology and are omitted for brevity.\\
\noindent{}\textbf{Step 2:} We use the same reasoning as above to linearize about $f_2$ to find $B_3$, which lies between the saddles $f_2$ and $f_3$:
\begin{align*}
     &B_2 = (\mathcal{O}(1),\mathcal{O}(1),\mathcal{O}(1) \delta^{(1-\varepsilon)\alpha-1},\mathcal{O}(1)\delta^{\alpha-1}) \rightarrow \\
     &B_3=(\mathcal{O}(1)\delta^{(\alpha-1)[(1-\varepsilon)\alpha-1]}, \mathcal{O}(1), \mathcal{O}(1),\mathcal{O}(1)\delta^{\varepsilon\alpha}).
\end{align*}
\textbf{Step 3:} We linearize about $f_3$ to find $B_4$, which lies between the saddles $f_3$ and $f_{34}$:
\begin{align*}
      &B_3=(\mathcal{O}(1)\delta^{(\alpha-1)[(1-\varepsilon)\alpha-1]}, \mathcal{O}(1), \mathcal{O}(1),\mathcal{O}(1)\delta^{\varepsilon\alpha}) \rightarrow \\
      &B_4 = (\mathcal{O}(1)\delta^{(\alpha-1)^2-\varepsilon\alpha^2},\mathcal{O}(1)\delta^{\varepsilon\alpha(\alpha-1)},\mathcal{O}(1),\mathcal{O}(1)).
\end{align*}
We require 
\begin{align*}
    \varepsilon<\varepsilon_{\scalebox{0.9}{$e$}} \text{, where } \varepsilon_{\scalebox{0.9}{$e$}} \coloneqq \frac{(\alpha-1)^2}{\alpha^2},
\end{align*}
so $u_4$ remains small. This requirement ensures (\ref{eqs:H4_step1_condition}) is satisfied.\\
\noindent{}\textbf{Step 4:} We linearize about $f_4$ to find $\hat{B}_1$, which lies between the saddles $f_{34}$ and $f_1$, thereby completing one revolution around the cycle:
\begin{align*}
    &B_4 = (\mathcal{O}(1)\delta^{(\alpha-1)^2-\varepsilon\alpha^2},\mathcal{O}(1)\delta^{\varepsilon\alpha(\alpha-1)},\mathcal{O}(1),\mathcal{O}(1)) \rightarrow \\
    &\hat{B}_1=(\mathcal{O}(1), \mathcal{O}(1)\delta^{\varepsilon\alpha(\alpha-1)+(1-\alpha+\varepsilon\alpha^2)[(\alpha-1)^2-\varepsilon\alpha^2]}, \mathcal{O}(1),\mathcal{O}(1)).
\end{align*}
We conclude that the heteroclinic cycle $H_{34}$ is stable if $\hat{v}_1<v_1$. It follows that the stability condition requires 
\begin{align*}
    &\varepsilon\alpha(\alpha-1)+(1-\alpha+\varepsilon\alpha^2)[(\alpha-1)^2-\varepsilon\alpha^2]>1.
\end{align*}
We find that the cycle is unstable for values of $\varepsilon$ between $\varepsilon_{**}$ and $\varepsilon_1$, where
\begin{align*}
    \varepsilon_1 \coloneqq 1 - \frac{2}{\alpha}.
\end{align*}
The ordering of these two critical values depends on $\alpha$. It is also important to remember that $f_{34}$ is only physical for 
\begin{align*}
    \varepsilon<\frac{1}{\alpha}=\varepsilon_{34}.
\end{align*}
Thus, the heteroclinic cycle $H_{34}$ ceases to exist for 
\begin{align*}
    \varepsilon>\text{min}(\varepsilon_{34},\varepsilon_{\scalebox{0.9}{$e$}}).
\end{align*}
\indent{}We reiterate that we must choose a sufficiently small coordinate $0<\delta\ll 1$ to obtain meaningful analytical results. The approximations of each successive point around the cycle rely on the time a trajectory spends near each saddle, which scales with $\ln{(1/\delta)}$. In the leading order approximation, we can safely neglect the contributions of the $\mathcal{O}(1)$ terms when $\delta$ is small enough that $\ln{(1/\delta)}$ dominates. Furthermore, this technique involves solving the linearized system both near and between saddles, in which $\mathcal{O}(1)$ terms become less reliable. We justify this method by observing $\delta \rightarrow 0$ as a trajectory approaches the $H_{34}$ cycle. Again, the contributions of $\mathcal{O}(1)$ terms can be neglected because the time spent near each saddle diverges to infinity, whereas the travel time---and corresponding error---between saddles becomes insignificant in comparison. This technique fails for the analysis of limit cycles because $\delta$ is not infinitesimal.\\

\noindent{}$\mathbf{H_{4}\textbf{: }f_1 \rightarrow f_2 \rightarrow f_3 \rightarrow f_{4} \rightarrow f_1.}$\\
\noindent{}\textbf{Step 1:} We follow the same procedure as for $H_{34}$. However, the initial point 
\begin{align*}
    B_1=(\mathcal{O}(1), \delta, \sigma,\mathcal{O}(1))
\end{align*}
lies between the saddles $f_{4}$ and $f_1$ and requires that $\delta$ and $\sigma$ both be small because the separatrix $f_{4} \rightarrow f_1$ lies in the subspace where $v=w=0$. By linearization of (\ref{eqs:3to4_system}) about $f_1$, we find 
\begin{align*}
    &B_1=(\mathcal{O}(1), \delta, \sigma,\mathcal{O}(1)) \rightarrow\\
    &B_2 = (\mathcal{O}(1),\mathcal{O}(1),\mathcal{O}(1) \sigma\delta^{(1-\varepsilon)\alpha-1},\mathcal{O}(1)\delta^{\alpha-1}).
\end{align*}
We require 
\begin{align*}
    \varepsilon<1-\frac{1}{\alpha} \text{ and } \alpha>1
\end{align*}
so $w_2$ and $z_2$ are small.\\
\noindent{}\textbf{Step 2:} As in $H_{34}$, we linearize about $f_2$ to find $B_3$, which lies between the saddles $f_2$ and $f_3$:
\begin{align*}
     &B_2 = (\mathcal{O}(1),\mathcal{O}(1),\mathcal{O}(1) \sigma\delta^{(1-\varepsilon)\alpha-1},\mathcal{O}(1)\delta^{\alpha-1}) \rightarrow \\
     &B_3=(\mathcal{O}(1)\sigma^{\alpha-1}\delta^{(\alpha-1)[(1-\varepsilon)\alpha-1]}, \mathcal{O}(1), \mathcal{O}(1),\mathcal{O}(1)\frac{1}{\sigma}\delta^{\varepsilon\alpha}).
\end{align*}
\textbf{Step 3:} We linearize about $f_3$ to find $B_4$, which lies between the saddles $f_3$ and $f_4$:
\begin{align*}
      &B_3=(\mathcal{O}(1)\sigma^{\alpha-1}\delta^{(\alpha-1)[(1-\varepsilon)\alpha-1]}, \mathcal{O}(1), \mathcal{O}(1),\mathcal{O}(1)\frac{1}{\sigma}\delta^{\varepsilon\alpha}) \rightarrow\\
      &B_4 = (\mathcal{O}(1)\sigma^\alpha\delta^{(\alpha-1)^2-\varepsilon\alpha^2},\mathcal{O}(1)\sigma^{-(\alpha-1)}\delta^{\varepsilon\alpha(\alpha-1)},\mathcal{O}(1),\mathcal{O}(1)).
\end{align*}
As with $H_{34}$, we require
\begin{align*}
    \varepsilon<\frac{(\alpha-1)^2}{\alpha^2}=\varepsilon_{\scalebox{0.9}{$e$}}
\end{align*}
so $u_4$ remains small. Again, this requirement ensures the condition imposed in Step 1 is satisfied.\\
\noindent{}\textbf{Step 4:} We linearize about $f_4$ to find $\hat{B}_1$, which lies between the saddles $f_{4}$ and $f_1$, thereby completing one revolution around the cycle:
\begin{align*}
    B_4 = ( &\mathcal{O}(1)\sigma^\alpha\delta^{(\alpha-1)^2-\varepsilon\alpha^2}, 
    \mathcal{O}(1)\sigma^{-(\alpha-1)}\delta^{\varepsilon\alpha(\alpha-1)}, \\
    &\mathcal{O}(1), \mathcal{O}(1) ) \rightarrow \\
    \hat{B}_1 = ( &\mathcal{O}(1), 
    \mathcal{O}(1)\sigma^{-(2\alpha-1)}\delta^{\varepsilon\alpha(2\alpha-1)-(\alpha-1)^2}, \\
    &\mathcal{O}(1)\sigma^{\alpha(\varepsilon\alpha-1)}\delta^{(\varepsilon\alpha-1)[(\alpha-1)^2-\varepsilon\alpha^2]}, 
    \mathcal{O}(1)).
\end{align*}
We first note that $\hat{w}_1$ is small when 
\begin{align*}
    \varepsilon>\frac{1}{\alpha}=\varepsilon_{34}.
\end{align*}
Therefore, $H_4$ only exists when 
\begin{align}\label{eqs:HH4_exists}
    \varepsilon_{34}<\varepsilon<\varepsilon_{\scalebox{0.9}{$e$}},
\end{align}
which is only possible when 
\begin{align*}
    \alpha>\frac{3+\sqrt{5}}{2}.
\end{align*}
Unlike $H_{34}$, where it was sufficient to compare the single small coordinate $\delta$ after one revolution, the stability condition for $H_4$ involves a relation between the two small coordinates $\delta$ and $\sigma$. We consider the logarithmic return map
\begin{align*}
    \vec{V}_1 &= 
    \begin{bmatrix}
        \ln{(\hat{v}_1)}\\
        \ln{(\hat{w}_1)}
    \end{bmatrix}\\[2ex]
    &=
    \renewcommand{\arraystretch}{2.0}
\setlength{\arraycolsep}{20pt}{
    \begin{bmatrix}
        \ln{(\sigma^{-(2\alpha-1)}\delta^{\varepsilon\alpha(2\alpha-1)-(\alpha-1)^2})}\\
        \ln{(\sigma^{\alpha(\varepsilon\alpha-1)}\delta^{(\varepsilon\alpha-1)[(\alpha-1)^2-\varepsilon\alpha^2]})}
    \end{bmatrix}}\\[2ex]
    & = 
    \renewcommand{\arraystretch}{2.0}
\setlength{\arraycolsep}{20pt}{
    \begin{bmatrix}
        [\varepsilon\alpha(2\alpha-1)-(\alpha-1)^2]\ln{(\delta)}+(1-2\alpha)\ln{(\sigma)}\\
        (\varepsilon\alpha-1)[(\alpha-1)^2-\varepsilon\alpha^2]\ln{(\delta)}+\alpha(\varepsilon\alpha-1)\ln{(\sigma)}
    \end{bmatrix}}\\[2ex]
    & =  
    \renewcommand{\arraystretch}{2.0}
\setlength{\arraycolsep}{20pt}{
    \begin{bmatrix}
        \varepsilon\alpha(2\alpha-1)-(\alpha-1)^2 &1-2\alpha\\
        (\varepsilon\alpha-1)[(\alpha-1)^2-\varepsilon\alpha^2] &\alpha(\varepsilon\alpha-1)
    \end{bmatrix}}
    \begin{bmatrix}
        \ln{(\delta)}\\
        \ln{(\sigma)}
    \end{bmatrix}\\[2ex]
    & = \mathbf{A}\vec{V}_0.
\end{align*}
Therefore, we find that 
\begin{align*}
    \vec{V}_n = \mathbf{A}\vec{V}_{n-1},
\end{align*}
for $n=1,2,3,\dots$, and 
\begin{align*}
    \vec{V}_n = \mathbf{A}^n\vec{V}_{0}.
\end{align*}
We chose $\delta$ and $\sigma$ to be small, which means the components of $\vec{V}_{0}$ are large negative numbers. We
conclude that the heteroclinic cycle is stable if the components of $\vec{V}_n$ tend to negative infinity as $n \rightarrow \infty$. To determine the necessary stability conditions, we consider a general matrix 
\begin{align*}
    \mathbf{A} = \begin{bmatrix}
        a & b\\
        c & d
    \end{bmatrix}
\end{align*}
with eigenvalues $\lambda_1$ and $\lambda_2$ and corresponding eigenvectors
\begin{align*}
    \begin{bmatrix}
        d-\lambda_1\\
        -c
    \end{bmatrix} \text{ and }     \begin{bmatrix}
        d-\lambda_2\\
        -c
    \end{bmatrix}.
\end{align*}
We assume both $\lambda_1$ and $\lambda_2$ are real, $\lambda_1>\lambda_2$, and refer the reader to Appendix \ref{apdx:H4} for the proof that $H_4$ is unstable when $\mathbf{A}$ has complex conjugate eigenvalues.\\
\indent{}To analyze $\mathbf{A}^n$, we diagonalize $\mathbf{A}$ to obtain
\begin{align*}
    &\mathbf{A}= P D P^{-1}=\\
    &\frac{1}{c(\lambda_1-\lambda_2)}
\begin{bmatrix} d-\lambda_1 & d-\lambda_2  \\ -c  & -c \end{bmatrix}
 \begin{bmatrix} \lambda_1 & 0 \\ 0  & \lambda_2 \end{bmatrix}
 \begin{bmatrix} -c & -(d-\lambda_2)\\ c & d-\lambda_1   \end{bmatrix}.
\end{align*}
Then,
\begin{align*}
&\mathbf{A}^n= P D^n P^{-1}= &&\\
&\frac{1}{c(\lambda_1-\lambda_2)}
\begin{bmatrix}
d-\lambda_1 & d-\lambda_2  \\
-c  & -c 
\end{bmatrix}
 \begin{bmatrix} 
 \lambda_1^n & 0 \\
 0  & \lambda_2^n 
 \end{bmatrix}
 \begin{bmatrix} 
 -c & -(d-\lambda_2)\\
 c & d-\lambda_1   \end{bmatrix}&&\\[2ex]
 &= \frac{1}{c(\lambda_1-\lambda_2)}\\[1ex]
 \renewcommand{\arraystretch}{2.2} 
&\begin{bmatrix} 
-c(d-\lambda_1)\lambda_1^n + c(d-\lambda_2)\lambda_2^n & 
-(d-\lambda_1)(d-\lambda_2)(\lambda_1^n-\lambda_2^n)\\[1ex]
c^2 (\lambda_1^n-\lambda_2^n)& c \lambda_1^n (d-\lambda_2)-c \lambda_2^n (d-\lambda_1)
\end{bmatrix},&&\\[1ex]
&\text{and }&&\\[1ex]
&\vec{V}_n = A^n \vec{V}_0=\frac{1}{c(\lambda_1-\lambda_2)}\begin{bmatrix} 
v_1 \\[1ex]
v_2
\end{bmatrix},&&
\end{align*}
where
\begin{align*}
    v_1 = &\ [-c(d-\lambda_1)\lambda_1^n + c(d-\lambda_2)\lambda_2^n]\ln (\delta) -\\
&(d-\lambda_1)(d-\lambda_2)(\lambda_1^n-\lambda_2^n) \ln (\sigma),\\
v_2 = &\ [c^2 (\lambda_1^n-\lambda_2^n)] \ln(\delta) +
 [c \lambda_1^n (d-\lambda_2)-c \lambda_2^n (d-\lambda_1)]\ln(\sigma).
\end{align*}
To have the components of $\vec{V}_n$ tend to negative infinity as $n \rightarrow \infty$, we require $\lambda_1>1$, such that terms grow in magnitude to infinity. The components of $\vec{V}_0$ are negative, so we must also require $|\lambda_2|<\lambda_1$, which guarantees that $\lambda_1$ dominates growth even if $\lambda_2<0$. It becomes sufficient to look only at the terms multiplied by $\lambda_1^n$ since it is the controlling term for large $n$. Thus,
\begin{align*}
    \vec{V}_n &\sim 
\frac{\lambda_1^n}{c(\lambda_1-\lambda_2)}
\renewcommand{\arraystretch}{2.0}
\setlength{\arraycolsep}{20pt}{
\begin{bmatrix} 
-c(d-\lambda_1) \ln (\delta) -
(d-\lambda_1)(d-\lambda_2)\ln (\sigma)\\
c^2 \ln(\delta) +
 c  (d-\lambda_2)\ln(\sigma)
\end{bmatrix}}\\[2ex]
&=
\frac{\lambda_1^n \left [ c\ln(\delta) +
   (d-\lambda_2)\ln(\sigma) \right ]}{\lambda_1-\lambda_2} 
   \renewcommand{\arraystretch}{1.5}
\setlength{\arraycolsep}{10pt}{
   \begin{bmatrix} 
-\frac{1}{c}(d-\lambda_1) \\
1
\end{bmatrix} }.
\end{align*}
Both components of $\vec{V}_n$ must go to negative infinity, so we require the components of the above vector to have the same sign. That is, we require
\begin{align}\label{eqs:lambda1}
\frac{1}{c}(d-\lambda_1)<0.
\end{align}
Consequently, it is necessary that
\begin{align}\label{eqs:lambda2}
    c\ln(\delta) +
   (d-\lambda_2)\ln(\sigma)<0.
\end{align}
In our model, when (\ref{eqs:HH4_exists}) is satisfied so that $H_4$ exists, we find that 
\begin{align*}
    c=(\varepsilon\alpha-1)[(\alpha-1)^2-\varepsilon\alpha^2]>0.
\end{align*}
It follows that we must have 
\begin{align*}
    \lambda_1>d=\alpha(\varepsilon\alpha-1)
\end{align*}
to satisfy (\ref{eqs:lambda1}). We can also see that (\ref{eqs:lambda2}) is unconditionally true for
\begin{align}\label{lambda2_less_than_d}
    \lambda_2-\alpha(\varepsilon\alpha-1)<0,
\end{align}
meaning that trajectories approach $H_4$ for any small $\delta$ and $\sigma$. In contrast, when
\begin{align*}
    \lambda_2-\alpha(\varepsilon\alpha-1)>0,
\end{align*}
$H_4$ is only attracting for trajectories with initial conditions that satisfy 
\begin{align}\label{eqs:HH4_ICs}
    \delta<\sigma^{\frac{\lambda_2-d}{c}}=\sigma^{\frac{\lambda_2-\alpha(\varepsilon\alpha-1)}{(\varepsilon\alpha-1)[(\alpha-1)^2-\varepsilon\alpha^2]}}.
\end{align}
In fact, when (\ref{lambda2_less_than_d}) holds, we find that (\ref{eqs:HH4_ICs}) is satisfied for all values of $\delta$ because the species population densities are restricted to values between zero and unity. Our analysis illustrates that a trajectory is only attracted to $H_4$ if $\varepsilon$ and $\alpha$ satisfy the required stability conditions \textit{and} the trajectory begins in the correct domain of phase space (i.e., the basin of attraction). \\
\indent{}We summarize the necessary conditions for the eigenvalues of $\mathbf{A}$, $\lambda_1$ and $\lambda_2$, under which $H_4$ is stable: (\MakeUppercase{\romannumeral 1}) the eigenvalues are both real (see Appendix \ref{apdx:H4}), (\MakeUppercase{\romannumeral 2}) $\lambda_1 > 1$, (\MakeUppercase{\romannumeral 3}) $\lambda_1 > \alpha(\varepsilon\alpha-1)$ and (\MakeUppercase{\romannumeral 4}) $|\lambda_2| < \lambda_1$. We note that the eigenvalues are always real in the $\alpha<3$ regime. In the $\alpha>3$ regime, we use $\varepsilon_{H_4}$ to indicate the transition from complex eigenvalues ($\varepsilon<\varepsilon_{H_4}$) to purely real eigenvalues ($\varepsilon>\varepsilon_{H_4}$), and the value of which is given by
\begin{align}\label{eqs:epsilon_HH}
    \varepsilon_{H_4} \coloneqq \frac{1}{\alpha(3\alpha-1)^2}(\beta+\omega),
\end{align}
where we define 
\begin{align*}
    &\beta = 5\alpha^3-10\alpha^2+10\alpha-3 \text{ and }\\
    &\omega = 2(\alpha-1)(2\alpha-1)\sqrt{(\alpha-1)(\alpha-3)}.
\end{align*}
Therefore, we find that $H_4$ is stable for 
\begin{align*}
    \varepsilon_{34}<\varepsilon<\varepsilon_{\scalebox{0.9}{$e$}} \text{ (for $\alpha<3$) and  } \varepsilon_{H_4}<\varepsilon<\varepsilon_{\scalebox{0.9}{$e$}} \text{ (for $\alpha>3$)},
\end{align*} 
when initial conditions satisfy (\ref{eqs:HH4_ICs}).\\

\noindent{}$\mathbf{H_{3}\textbf{: }f_1 \rightarrow f_2 \rightarrow f_3  \rightarrow f_1}$.\\
We note that the $H_{34}$ and $H_4$ cycles both require 
\begin{align*}
    \varepsilon<\frac{(\alpha-1)^2}{\alpha^2}=\varepsilon_{\scalebox{0.9}{$e$}}
\end{align*}
for $u_4$ to be small and for the trajectory to travel from $f_3$ to $f_{34}$ or $f_4$. If this condition is not met, then $u_4$ grows and $f_3$ is connected to $f_1$ (forming the cycle $H_3$).\\
\noindent{}\textbf{Step 1:} Again, we follow the same procedure as for $H_{34}$. This time, the initial point 
\begin{align*}
    B_1=(\mathcal{O}(1), \delta, \mathcal{O}(1),\sigma)
\end{align*}
lies between the saddles $f_{3}$ and $f_1$ and requires that $\delta$ and $\sigma$ both be small because the separatrix $f_{3} \rightarrow f_1$ lies in the subspace where $v=z=0$. By linearization of (\ref{eqs:3to4_system}) about $f_1$, we find 
\begin{align*}
    &B_1=(\mathcal{O}(1), \delta,\mathcal{O}(1),\sigma) \rightarrow \\
    &B_2 = (\mathcal{O}(1),\mathcal{O}(1),\mathcal{O}(1) \delta^{(1-\varepsilon)\alpha-1},\mathcal{O}(1)\sigma\delta^{\alpha-1}).
\end{align*}
We require 
\begin{align*}
    \varepsilon<1-\frac{1}{\alpha} \text{ and }\alpha>1
\end{align*}
so $w_2$ and $z_2$ are small.\\
\noindent{}\textbf{Step 2:} We linearize about $f_2$ to find $B_3$, which lies between the saddles $f_2$ and $f_3$:
\begin{align*}
     &B_2 = (\mathcal{O}(1),\mathcal{O}(1),\mathcal{O}(1) \delta^{(1-\varepsilon)\alpha-1},\mathcal{O}(1)\sigma\delta^{\alpha-1}) \rightarrow\\
     &B_3=(\mathcal{O}(1)\delta^{(\alpha-1)[(1-\varepsilon)\alpha-1]}, \mathcal{O}(1), \mathcal{O}(1),\mathcal{O}(1)\sigma\delta^{\varepsilon\alpha}).
\end{align*}
\noindent{}\textbf{Step 3:} We linearize about $f_3$ to find $\hat{B}_1$, which lies between the saddles $f_{3}$ and $f_1$, thereby completing one revolution around the cycle:
\begin{align*}
    B_3=(&\mathcal{O}(1)\delta^{(\alpha-1)[(1-\varepsilon)\alpha-1]}, \mathcal{O}(1), \mathcal{O}(1),\mathcal{O}(1)\sigma\delta^{\varepsilon\alpha}) \rightarrow \\
    \hat{B}_1=(&\mathcal{O}(1), \mathcal{O}(1)\delta^{(\alpha-1)^2[(1-\varepsilon)\alpha-1]},\\ &\mathcal{O}(1),\mathcal{O}(1)\sigma\delta^{\varepsilon\alpha-(\alpha-1)[(1-\varepsilon)\alpha-1]}).
\end{align*}
We conclude that the cycle is stable when $\hat{v}_1<v_1$ and $\hat{z}_1<z_1$. Thus, we require
\begin{align*}
    (\alpha-1)^2[(1-\varepsilon)\alpha-1]>1
\end{align*}
and
\begin{align*}
    \varepsilon\alpha-(\alpha-1)[(1-\varepsilon)\alpha-1]>0.
\end{align*}
It follows that $H_3$ is stable for 
\begin{align*}
    &\varepsilon>\frac{(\alpha-1)^2}{\alpha^2}=\varepsilon_{\scalebox{0.9}{$e$}} \text{ (as expected)}\\
    &\text{and } \varepsilon < \frac{(\alpha-1)^3-1}{\alpha(\alpha-1)^2}=\varepsilon_{123}.
\end{align*}

\begin{figure}[]
    \centering
    \includegraphics[width=1\linewidth]{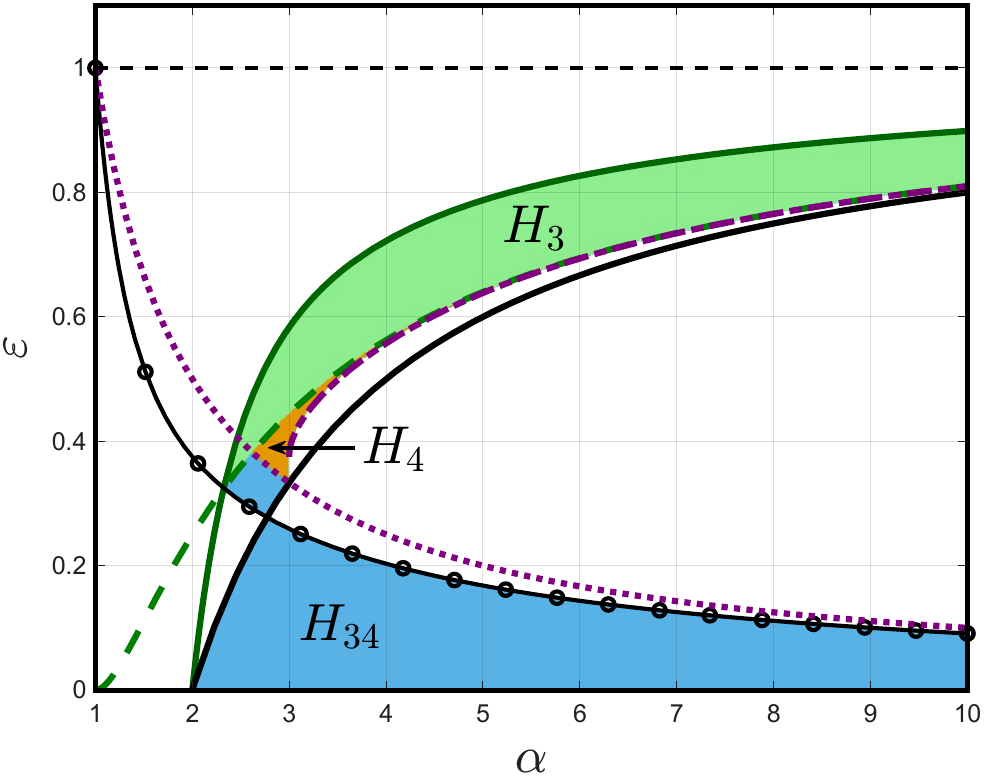}
    \caption{Regions of heteroclinic cycle stability in the $\alpha\text{-}\varepsilon$ plane, where $\alpha$ represents inter-species competition and $\varepsilon$ is the homotopy parameter. The colored regions indicate the following cycles are stable: $H_{34}$ (blue), $H_4$ (orange), and $H_3$ (green). The critical $\varepsilon$ values are as follows: $\varepsilon_{**}$ (solid line with circles; where the system obtains the curve of steady states $f_{1234a}$), $\varepsilon_{1}$ (solid black line), $\varepsilon_{123}$ (solid green line), $\varepsilon_{\scalebox{0.9}{$e$}}$ (dashed green line), $\varepsilon_{34}$ (dotted purple line), and $\varepsilon_{H_4}$ (dashed-dotted purple line).}
    \label{fig:Regions_of_HCycle_Stability_3to4Model}
\end{figure}

\noindent{}\textbf{Summary.} As in the RPS model, there are no stable heteroclinic cycles when $\alpha<2$. Moreover, when $\alpha>2$ and $\varepsilon$ is small, our model remains structurally similar to an RPS model and can exhibit RPS-like cycles of competitive exclusion. We note that in our model, the absence of species $u$ allows both species $z$ and $w$ to thrive, which is why we observe cycles involving four species, rather than $H_3$, when $\varepsilon$ is small. Furthermore, small values of $\varepsilon$ mean that the competitive pressure exerted by $z$ on $w$ is low and the two species can coexist through the saddle $f_{34}$, resulting in the stable $H_{34}$ cycle.\\
\indent{}As inter-species competition increases and $\alpha \rightarrow \infty$, the stability regions for the heteroclinic cycles contract and eventually vanish (Fig. \ref{fig:Regions_of_HCycle_Stability_3to4Model}). For example, $H_{34}$ exists and is stable for $\varepsilon<\text{min}(\varepsilon_{34}, \varepsilon_{\scalebox{0.9}{$e$}})$, but $\varepsilon_{34}\rightarrow 0$ as $\alpha\rightarrow \infty$. Thus, $H_{34}$ is only stable for values of $\varepsilon$ increasingly closer to zero, while $H_4$ and $H_3$ are stable for values of $\varepsilon$ tending towards unity. We refer the reader to Table \ref{tab:critical_epsilon_table} for a complete list of critical $\varepsilon$ values and to Table \ref{tab:summary_of_states_table} for a summary of all stability regions.\\
\indent{}We find that it becomes impossible for all four species to coexist when inter-species competition is high (Fig. \ref{fig:Regions_of_HCycle_Stability_3to4Model}). Instead, we see two possible behaviors depending on the strength of species interactions (i.e., the value of $\varepsilon$). When values of $\varepsilon$ are approximately zero or unity and our modified model behaves like the RPS or Hand of Bridge systems, species must self-segregate in time via heteroclinic cycles or by forming noncompetitive alliances. Consequently, we see cycling behavior. In contrast, for intermediate values of $\varepsilon$, the combined pressures from species $u$ and $z$ are too high for species $w$ to survive. The weakened $w$ population allows $v$ to grow and form the stable $v\text{-}z$ alliance.\\
\indent{}We note that a narrow window of limit cycles emerges for intermediate values of $\alpha$, which we discuss in more detail below. 

\renewcommand{\arraystretch}{1.8}
\begin{table}[]
\centering
\caption{Summary table of critical $\varepsilon$ values.}
\label{tab:critical_epsilon_table}
\begin{tabular*}{\tblwidth}{@{}L C p{4.5cm}@{}}
    \toprule
    \textbf{Name } & \textbf{Formula} & \textbf{Meaning} \\ 
    \midrule
    $\varepsilon_*$ & $1-\dfrac{4}{\alpha}+\dfrac{8}{\alpha^2}-\dfrac{8}{\alpha^3}$ \ & Stability boundary of $f_{1234s}$ \\ 
    \addlinespace
    $\varepsilon_{**}$ & $\dfrac{1}{\alpha}-\dfrac{1}{\alpha^2}+\dfrac{1}{\alpha^3}$ & Stability boundary of $f_{1234s}$, $f_{123}$, $f_{234}$, and $H_{34}$ \\ 
    \addlinespace
    $\varepsilon_{123}$ & $\dfrac{(\alpha-1)^3-1}{\alpha(\alpha-1)^2}$ & Stability boundary of $f_{123}$ and $H_3$ \\ 
    \addlinespace
    $\varepsilon_{13}$ & $1-\dfrac{1}{\alpha}+\dfrac{1}{\alpha^2}$ & Stability boundary of $f_{123}$ and $f_{13}$ \\ 
    \addlinespace
    $\varepsilon_{34}$ & $\dfrac{1}{\alpha}$ & \vspace{-0.5\baselineskip}\begin{itemize}[nosep, leftmargin=*, label=-]
        \item Stability boundary of $f_{234}$ and $f_{24}$
        \item The state $f_{34}$ becomes unphysical for $\varepsilon>\varepsilon_{34}$
        \item Stability transitions from $H_{34}$ to $H_4$
    \end{itemize}\\ 
    \addlinespace
    $\varepsilon_{1}$ & $1-\dfrac{2}{\alpha}$ & Stability boundary of $H_{34}$ \\ 
    \addlinespace
    $\varepsilon_{\scalebox{0.9}{$e$}}$ & $\dfrac{(\alpha-1)^2}{\alpha^2}$ & Stability transition from $H_{34}$ or $H_4$ to $H_3$ \\ 
    \addlinespace
    $\varepsilon_{H_4}$ & \text{See equation (\ref{eqs:epsilon_HH})} & Stability boundary of $H_4$ for $\alpha>3$ \\ 
    \bottomrule
\end{tabular*}
\end{table}

\renewcommand{\arraystretch}{1.8}
\begin{table}[]
\centering
\caption{Summary table of steady states and heteroclinic cycles. The state $f_{1234a}$ is not listed as this family exists only for $\varepsilon_{**}$.}
\label{tab:summary_of_states_table}
\begin{tabular}{@{}l l@{}}
    \toprule
    \textbf{Name} & \textbf{Stability Window}\\ 
    \midrule
    $f_{1234s}$ & $\varepsilon_* < \varepsilon<\varepsilon_{**}$\\
    \addlinespace
    $f_{123}$ & $\text{max}(\varepsilon_{123},\varepsilon_{**}) < \varepsilon<\varepsilon_{13}$\\
    \addlinespace
    $f_{234}$ & $\varepsilon_{**} < \varepsilon<\varepsilon_{34}$\\
    \addlinespace
    $f_{13}$ & $\varepsilon>\varepsilon_{13}$\\
    \addlinespace
    $f_{24}$ & $\varepsilon>\varepsilon_{34}$\\
    \addlinespace
    $H_{34}$ & 
    \begin{varwidth}[t]{\linewidth}$\displaystyle\varepsilon<\min(\varepsilon_{34},\varepsilon_{\scalebox{0.9}{$e$}})$
    \begin{itemize}[nosep, leftmargin=*, label=-] 
        \item $\varepsilon \notin (\varepsilon_1,\varepsilon_{**})$ for $\alpha<2.7693$
        \item $\varepsilon \notin (\varepsilon_{**},\varepsilon_1)$ for $\alpha>2.7693$
    \end{itemize}\end{varwidth}\\
    \addlinespace
    $H_4$ & \begin{varwidth}[t]{\linewidth}
    \begin{itemize}[nosep, leftmargin=*, label=-] 
        \item $\varepsilon_{34} < \varepsilon < \varepsilon_{\scalebox{0.9}{$e$}} $ for $\alpha<3$
        \item $\varepsilon_{H_4} < \varepsilon < \varepsilon_{\scalebox{0.9}{$e$}} $ for $\alpha>3$
    \end{itemize}\end{varwidth}\\
    \addlinespace
    $H_3$ & $\varepsilon_{\scalebox{0.9}{$e$}}<\varepsilon<\varepsilon_{123}$ \\
    \bottomrule
\end{tabular}
\end{table}

\section{Numerical Simulations}
\indent{}\indent{}We describe numerical simulations that confirm the analytical results presented above and demonstrate the existence of a window of limit cycles. We have identified the following regimes: the absence of limit cycles ($\alpha<2.7693$), a narrow window of limit cycles ($2.7693 < \alpha < 3$), and no stable coexistence ($\alpha>3$). We also discuss system behavior for large inter-species competition ($\alpha \rightarrow \infty$).\\
\indent{}We implement standard ODE solvers to obtain our numerical results and improve the accuracy of our simulations by introducing the following log-transform change of variables:
\begin{align*}
    u_l=\ln{u}, v_l=\ln{v}, w_l=\ln{w}, \text{ and } z_l=\ln{z}.
\end{align*}
Reformulating the system (\ref{eqs:3to4_system}) in this way allows us to track a solution trajectory for longer times as it gets progressively closer to the saddles comprising a heteroclinic cycle. Without the change of variables, a solution visiting a saddle becomes trapped at the fixed point when the small population densities fall below machine epsilon.\\
\indent{}To characterize the different heteroclinic cycles in each of the regimes, we consider the duration of visits by the solution trajectory to one of the saddle points. All heteroclinic cycles observed in this system involve a visit to the saddle $f_1$, so we implement the technique described in \cite{siegfried2025fractured} to track the number and duration of each visit by species $u$ to the $f_1$ saddle. We consider when $u$ approaches its carrying capacity (unity) from below and define a visit to the saddle when the population size of $u$ surpasses a prescribed parameter ($u_{\text{critical}}$) close to the carrying capacity (e.g., when $u>u_{\text{critical}}=0.9999)$. The duration of the $n$th visit to the saddle point is $T_n$, where the duration is defined to be the length of time for which $u>u_{\text{critical}}$. 
We use the scaling
\begin{align}\label{eqs:mu}
    T_n \approx Ae^{\mu n} \Rightarrow \ln{(T_n)} \approx \mu n + \ln(A),
\end{align}
where $A$ depends on the choice of saddle and the definition of $u_{\text{critical}}$. When $n$ is large, $T_n$ is well approximated for any saddle choice, since the durations of the visits to any saddle connected by a heteroclinic cycle expand at the same rate as time tends to infinity. We determine $\mu$ by applying a linear least-squares fit to the log-transformed visit durations, using the scaling defined in (\ref{eqs:mu}). Thus, we use the value $\mu$ to characterize the heteroclinic cycles. That is, positive values of $\mu$ indicate that a cycle is attracting. Additionally, strongly attracting (or ``sticky'') heteroclinic cycles are characterized by large values of $\mu$, while $\mu$ is small for weakly attracting cycles. We note that $\mu=0$ (within machine tolerance) indicates that the cycling behavior is due to a periodic limit cycle, rather than a weakly attracting heteroclinic cycle.\\

\begin{figure}[]
    \centering
    \includegraphics[width=1\linewidth]{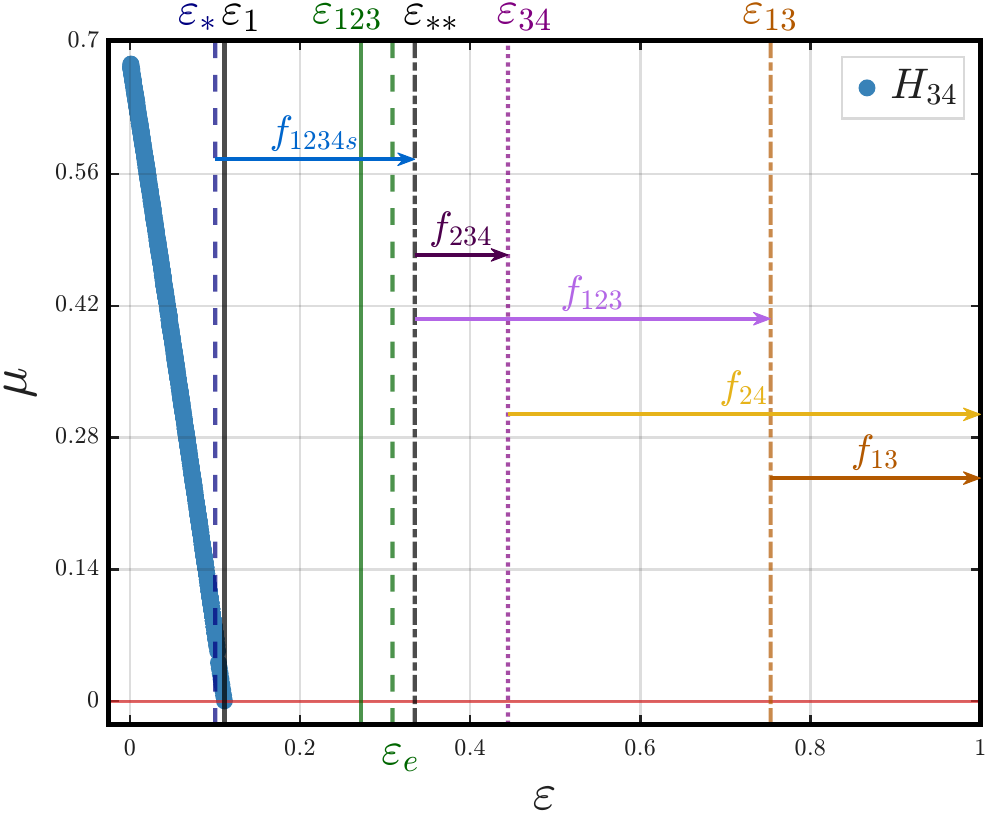}
    \caption{Heteroclinic cycle growth rate $\mu$ as a function of the homotopy parameter $\varepsilon$ for $\alpha=2.25$, characterizing $H_{34}$ (blue). The vertical lines indicate important values of $\varepsilon$, and horizontal arrows show regions where steady states are stable. Dynamic behavior ends when the blue circles hit zero and $H_{34}$ loses stability.}
    \label{fig:ep_vs_mu_2_25}
\end{figure}

\begin{figure}[]
    \centering
    \includegraphics[width=1\linewidth]{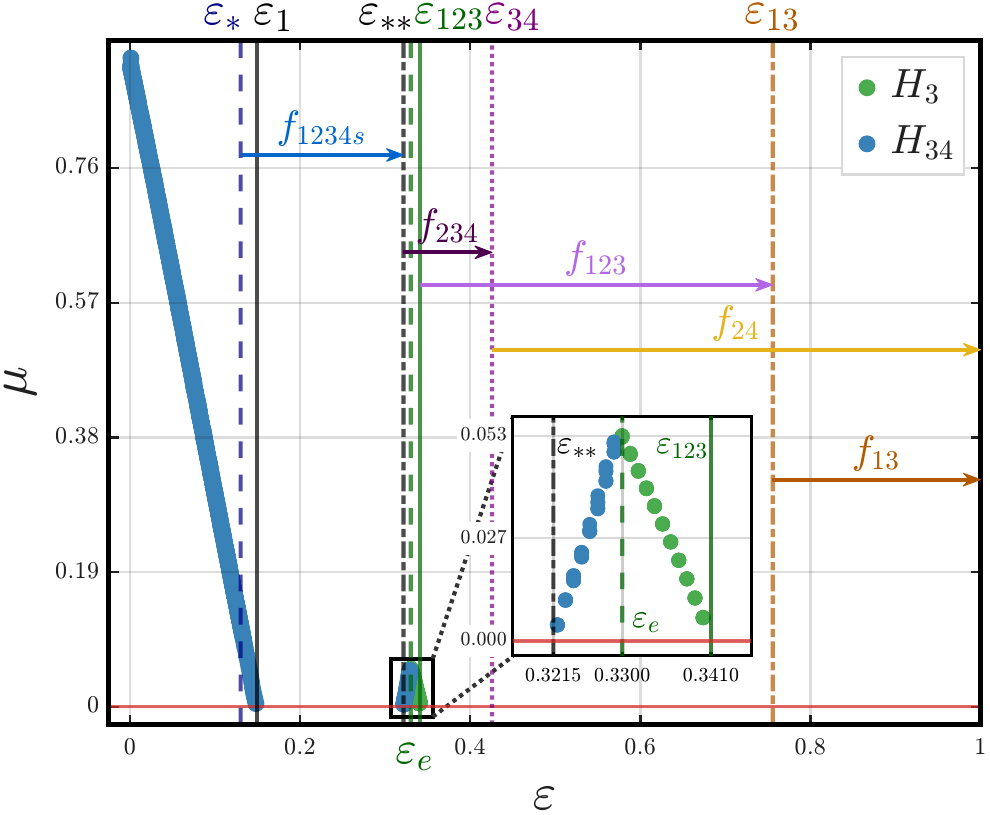}
    \caption{Heteroclinic cycle growth rate $\mu$ as a function of the homotopy parameter $\varepsilon$ for $\alpha=2.35$, characterizing $H_{34}$ (blue) and $H_3$ (green). The vertical lines indicate important values of $\varepsilon$, and horizontal arrows show regions where steady states are stable.}
    \label{fig:ep_vs_mu_2_35}
\end{figure}

\textbf{No limit cycles ($\boldsymbol{\alpha<}\mathbf{2.7693}$)}. Within this regime, we do not observe any limit cycles. However, we note that there are three distinct sub-regimes, each pertaining to the behavior of heteroclinic cycles.\\
\indent{}First, numerical simulations show agreement with our analytical results that there are no heteroclinic cycles when $\alpha<2$ (\MakeUppercase{\romannumeral 1}). In the sub-regime $2<\alpha<2.3247$ (\MakeUppercase{\romannumeral 2}), the only stable cycle is $H_{34}$. Calculations of $\mu$ show that $H_{34}$ is most strongly attracting when $\varepsilon=0$, and the strength of this attraction monotonically decreases until the cycle loses stability at $\varepsilon_1$ (Fig. \ref{fig:ep_vs_mu_2_25}). As $\varepsilon$ increases to unity, the system no longer exhibits cycling behavior and, instead, transitions towards bistability between alliances of non-competing species---fundamental Hand of Bridge behavior. In the sub-regime $2.3247 < \alpha < 2.7693$ (\MakeUppercase{\romannumeral 3}), the transition from the coexistence steady state to three- and then two-species equilibria is interrupted by the reappearance of stable heteroclinic cycles as $\varepsilon$ increases to unity. In this sub-regime, the coexistence steady state becomes an unstable focus and $H_{34}$ regains stability at $\varepsilon_{**}$. Solution trajectories are increasingly attracted to the cycle until a local maximum is reached at $\varepsilon_{\scalebox{0.9}{$e$}}$, at which point stability is transferred to $H_3$ (Fig. \ref{fig:ep_vs_mu_2_35}). At $\varepsilon_{123}$, all cycles become unstable, and the system transitions towards bistability between $f_{13}$ and $f_{24}$. That is, all dynamic behavior is quenched on the route to $\varepsilon=1$, resulting in the dynamically stable Hand of Bridge system\\

\begin{figure}[]
    \centering
    \includegraphics[width=1\linewidth]{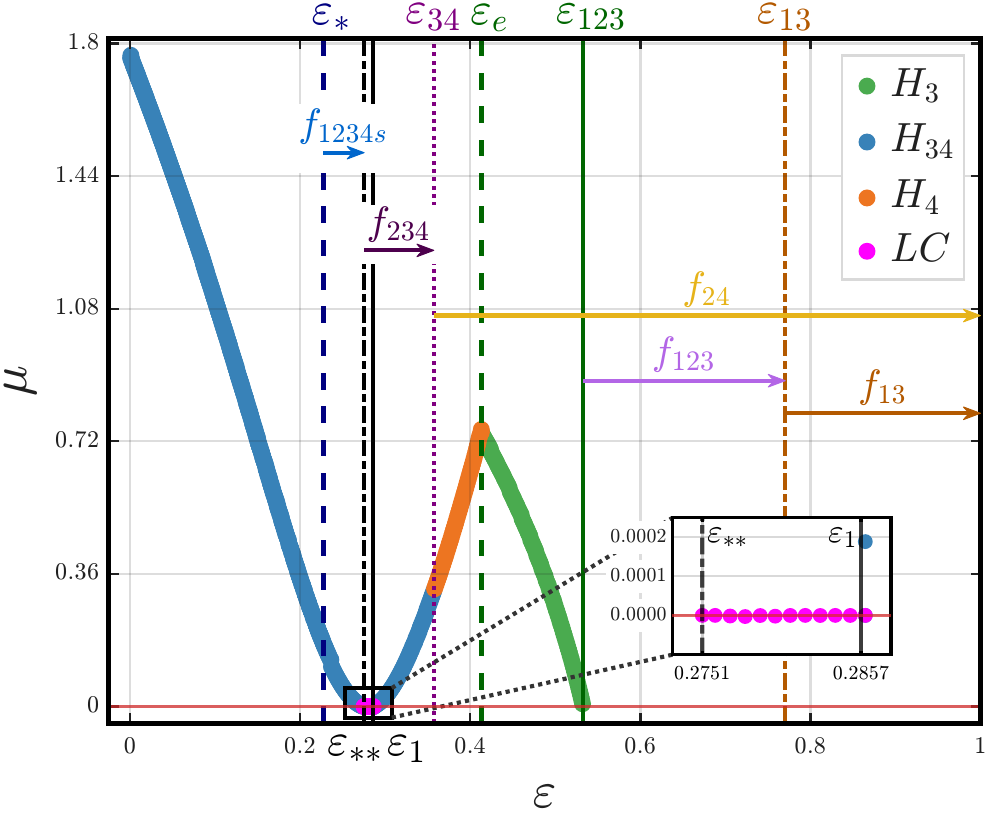}
    \caption{Heteroclinic cycle growth rate $\mu$ as a function of the homotopy parameter $\varepsilon$ for $\alpha=2.80$, characterizing $H_{34}$ (blue), limit cycles ($LC$, magenta), $H_4$ (orange), and $H_3$ (green). The vertical lines indicate important values of $\varepsilon$, and horizontal arrows show regions where steady states are stable.}
    \label{fig:ep_vs_mu_2_8}
\end{figure}

\textbf{A window of limit cycles ($\mathbf{2.7693} \boldsymbol{<\alpha<}\mathbf{3}$)}. At $\alpha=2.7693$, the critical values $\varepsilon_1$ and $\varepsilon_{**}$ intersect. Beyond this intersection ($\alpha>2.7693$), the order of the two values changes from $\varepsilon_1<\varepsilon_{**}$ to $\varepsilon_{**}<\varepsilon_1$. In the previous regime, $f_{1234s}$ was stable while $H_{34}$ was unstable for values of $\varepsilon \in (\varepsilon_1,\varepsilon_{**})$. In this new regime, both $f_{1234s}$ and $H_{34}$ are unstable in the parameter window $(\varepsilon_{**},\varepsilon_1)$. Specifically, it is in this regime that we see a local steady-state bifurcation around $f_{1234s}$, the creation of a line of steady states $f_{1234a}$, and a global heteroclinic bifurcation all at the same critical value $\varepsilon_{**}$. We see a family of periodic orbits at $\varepsilon_{**}$ and the existence of a limit cycle for a small range of $\varepsilon>\varepsilon_{**}$ (Fig. \ref{fig:ep_vs_mu_2_8}).

\begin{figure}[]
    \centering
    \begin{subfigure}[t]{0.45\textwidth}
        \centering
        \includegraphics[width=\textwidth]{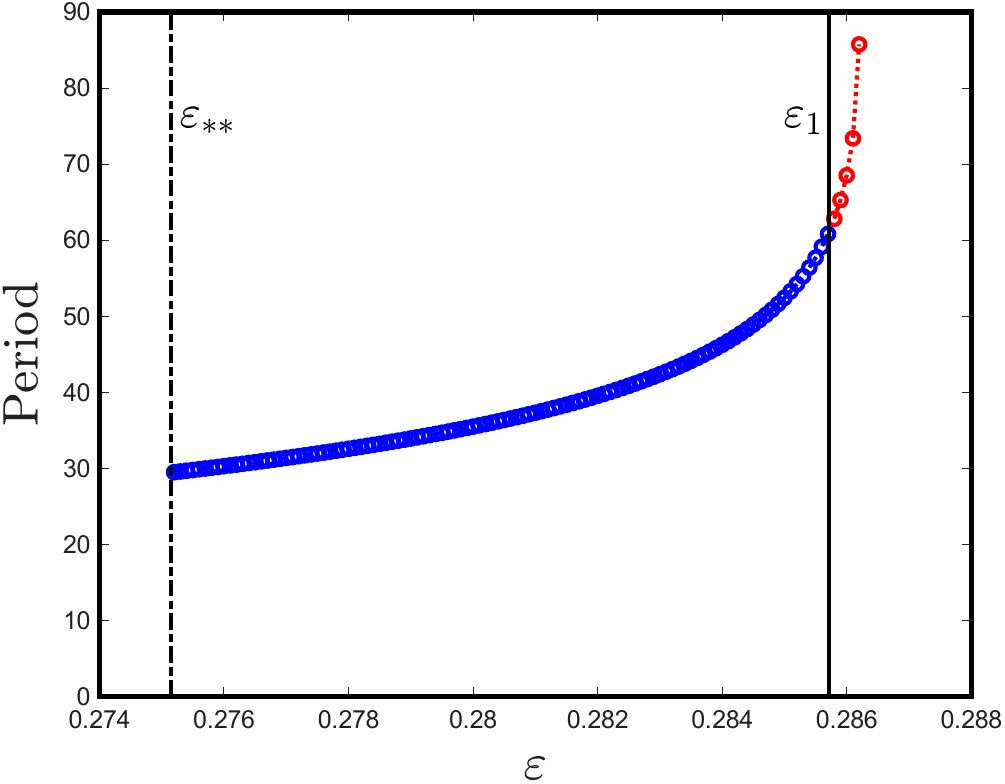}
        \caption{Periods for the observed limit cycles. Blue dots represent limit cycles observed within the parameter window where $H_{34}$ is unstable, while red dots represent limit cycles observed when $H_{34}$ regains stability.\\}
        \label{subfig:LC_Periods}
    \end{subfigure}
    \vfill
    \begin{subfigure}[t]{0.45\textwidth}
        \centering
        \includegraphics[width=\textwidth]{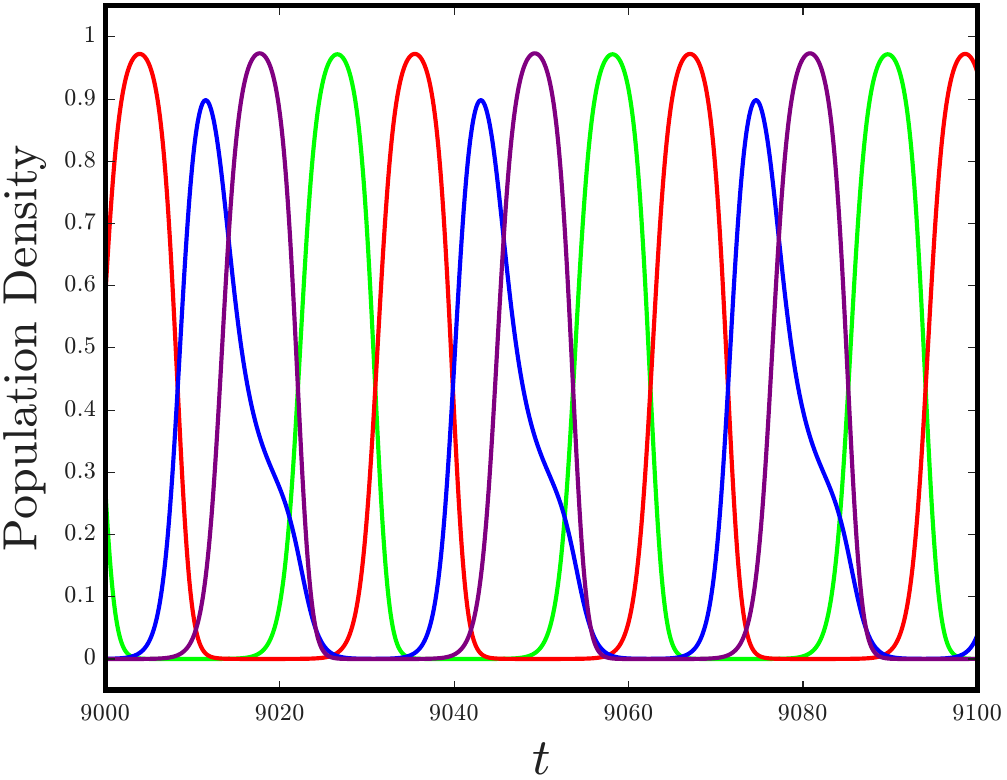}
        \caption{Limit cycle solution for $\varepsilon=0.27516$, showing population density for species $u$ (green), $v$ (red), $w$ (blue) and $z$ (purple).}
        \label{subfig:LC_example}
    \end{subfigure}
    \vfill
    \caption{Limit cycles observed when $\alpha=2.8$.}
    \label{fig:LC_Periods_and_example}
\end{figure}

\begin{figure}[]
    \centering
    \begin{subfigure}[t]{0.49\textwidth}
        \centering
        \includegraphics[width=\textwidth]{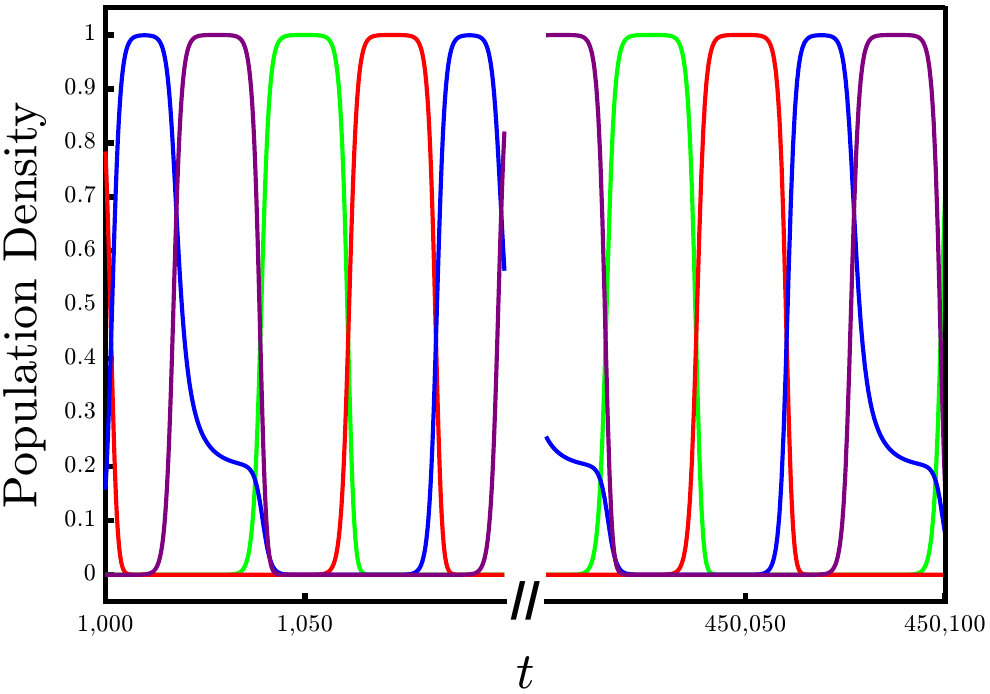}
        \caption{Limit cycle solution for $\varepsilon=0.2862$, showing population density for species $u$ (green), $v$ (red), $w$ (blue) and $z$ (purple).\\}
        \label{subfig:LC_ep2862}
    \end{subfigure}
    \vfill
    \begin{subfigure}[t]{0.49\textwidth}
        \centering
        \includegraphics[width=\textwidth]{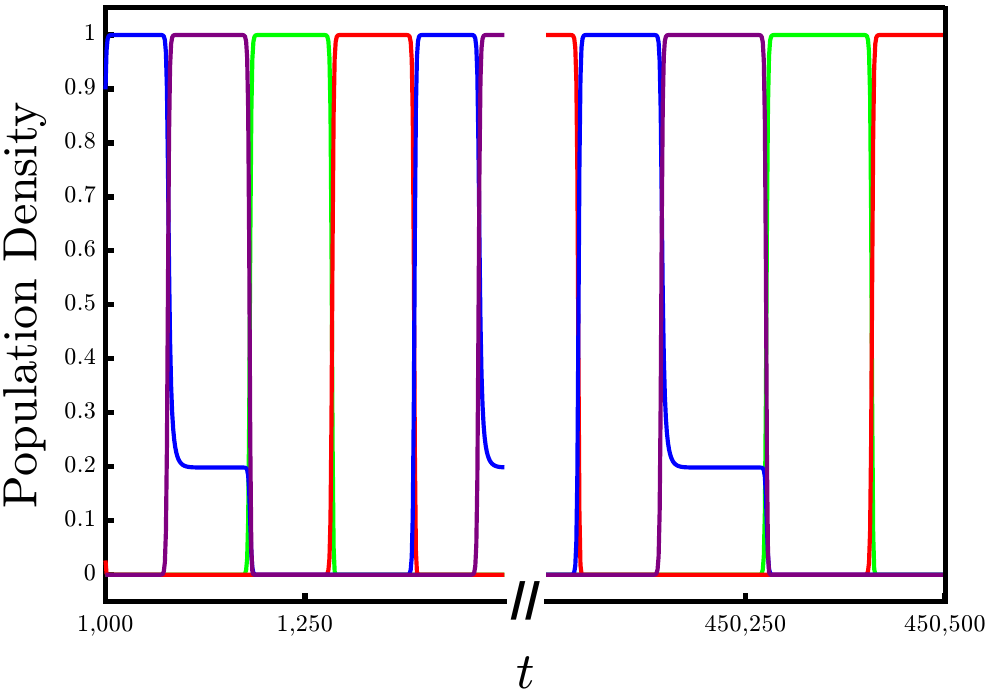}
        \caption{Heteroclinic cycle $H_{34}$ for $\varepsilon=0.2862$, showing population density for species $u$ (green), $v$ (red), $w$ (blue) and $z$ (purple).}
        \label{subfig:H4_ep2862}
    \end{subfigure}
    \vfill
    \caption{Limit cycle and $H_{34}$ coexist for $\varepsilon=0.2862$ and $\alpha=2.8$. Observe the different time scales between the two $t$-axes: plot (a) increments by $50$ units of $t$ while plot (b) increments by $250$ units.}
    \label{fig:LC_and_H4_ep2862}
\end{figure}

\indent{}Figure \ref{subfig:LC_Periods} plots the calculated periods of the observed limit cycles against values of $\varepsilon$ (not including the family of limit cycles at $\varepsilon_{**}$), and Figure \ref{subfig:LC_example} shows the limit cycle solution for $\varepsilon=0.27516$. Moreover, we observed limit cycles for values of $\varepsilon$ greater than $\varepsilon_1$, indicating that there is a region of tri-stability between a limit cycle, $H_{34}$, and $f_{234}$ (Fig. \ref{subfig:LC_Periods}). We note that the limit cycle solutions appear visually distinct from the $H_{34}$ solutions (Fig. \ref{fig:LC_and_H4_ep2862}).\\
\indent{}The cycle $H_{34}$ is very weakly attracting for values of $\varepsilon$ near $\varepsilon_1$. To compute the $H_{34}$ cycle in the region of tri-stability, we first computed the limit cycle solution for a given $\varepsilon>\varepsilon_1$ and determined the minimum values of $u,w$ and $z$ (called $u_\text{min}, w_\text{min}$, and $z_\text{min}$). We then took the initial conditions to be
\begin{align*}
    (u_\text{min},V, w_\text{min},z_\text{min}),
\end{align*}
where $V$ is several orders of magnitude smaller than the other three minimum values. As $\varepsilon$ increases, solution trajectories become more strongly attracted to $H_{34}$ and we no longer observe limit cycles (either because the limit cycle becomes unstable or because the basin of attraction shrinks to a region too small for us to numerically identify). The largest value of $\varepsilon$ for which we identified the existence of a limit cycle for $\alpha = 2.8$ is $\varepsilon=0.2862$ (Fig. \ref{subfig:LC_ep2862}).\\

\begin{figure}[]
    \centering
    \includegraphics[width=1\linewidth]{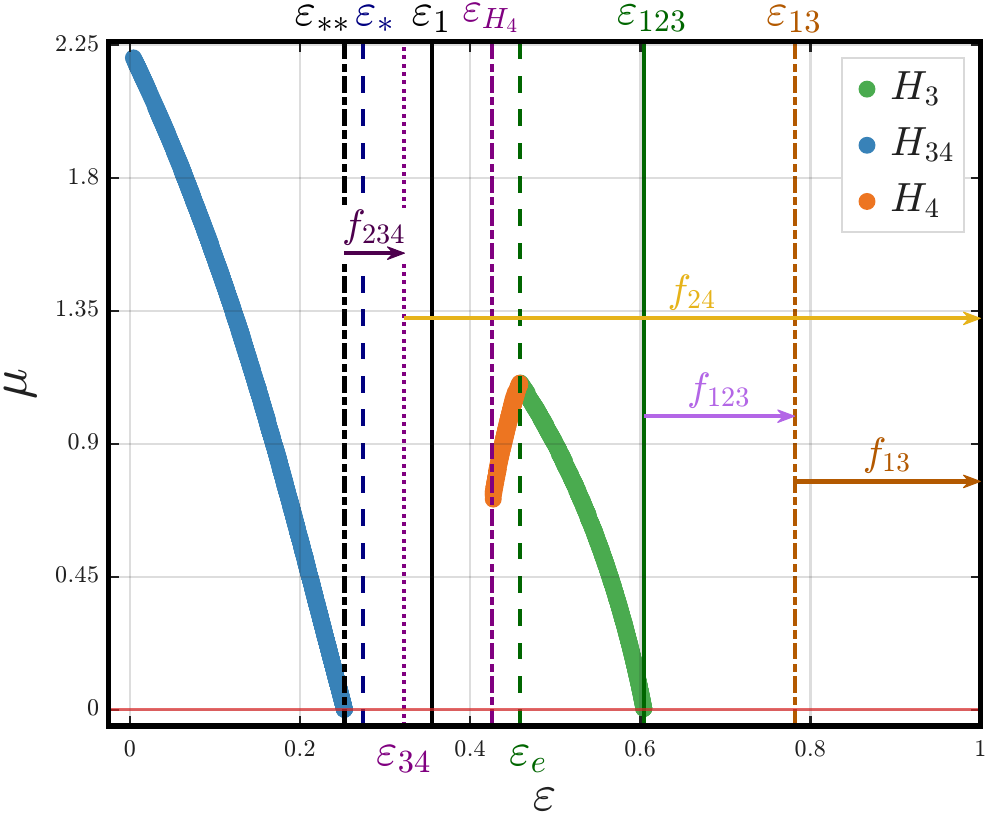}
    \caption{Heteroclinic cycle growth rate $\mu$ as a function of the homotopy parameter $\varepsilon$ for $\alpha=3.10$, characterizing $H_{34}$ (blue), $H_4$ (orange), and $H_3$ (green). The vertical lines indicate important values of $\varepsilon$, and horizontal arrows show regions where steady states are stable.}
    \label{fig:ep_vs_mu_3_1}
\end{figure}

\textbf{No stable coexistence ($\boldsymbol{\alpha>}\mathbf{3}$).} At the point $\alpha=3$, the curves $\varepsilon_*$ and $\varepsilon_{**}$ intersect at $(\alpha,\varepsilon)=(3,0.2593)$, as well as the curves $\varepsilon_1$ and $\varepsilon_{34}$ at $(\alpha,\varepsilon)=(3,\frac{1}{3})$. The system exhibits two distinct codimension-2 bifurcation points, and is highly sensitive to the initial conditions, the choice of numerical solver, and numerical tolerances. Beyond this value of $\alpha$, in the regime $\alpha>3$, the coexistence steady state $f_{1234s}$ is always unstable. It is also in this regime that we observe stable steady states only ($f_{234}$ and $f_{24}$) rather than the coexistence of $f_{234}$ and limit cycles for values of $\varepsilon \in (\varepsilon_{**},\varepsilon_1)$ (Fig. \ref{fig:ep_vs_mu_3_1}). That is, the limit cycles that emerged in the previous regime are no longer observed for $\alpha>3$.\\ 
\indent{}We reiterate that as inter-species competition increases ($\alpha \rightarrow \infty$), cyclically interacting species struggle to coexist and must, instead, self segregate in time via cycles or alliances. As discussed below, we find that this regime favors the alliance state $f_{24}$. \\

\textbf{Strong Inter-Species Competition ($\boldsymbol{\alpha \rightarrow \infty}$).} When inter-species competition is strong (i.e., $\alpha$ is large), species $w$ is weakened and vulnerable to extinction for intermediate values of $\varepsilon$ due to the combined pressure from species $u$ and $z$. Since the heteroclinic cycles $H_{34}$, $H_4$, and $H_3$ rely on the survival of species $w$, large values of $\alpha$ push the stability regions away from intermediate $\varepsilon$. That is, $w$ cannot survive for intermediate $\varepsilon$ and, consequently, we only see stable heteroclinic cycles for $\varepsilon$ near zero or unity. Similarly, the stable states $f_{123}$, $f_{234}$, and $f_{13}$ necessitate the survival of species $w$. As the survival of $w$ becomes more precarious with increasing $\alpha$, the states involving $w$ require more extreme values of $\varepsilon$ to remain stable. \\
\indent{}For example, as $\alpha \rightarrow \infty$, it becomes difficult to sustain $H_4$, $H_3$, $f_{123}$, and $f_{13}$ at intermediate values of $\varepsilon$. The presence of species $u$ in these configurations restricts their stability to the extreme $\varepsilon \sim 1$ boundary, at which $w$ loses its vulnerability to $u$. Likewise, $H_{34}$ and $f_{234}$ are only stable when $\varepsilon \sim 0$ because species $z$ and $w$ can coexist only when $w$ loses its vulnerability to $z$. Thus, strong inter-species competition drives a preference for the stable state $f_{24}$ at intermediate $\varepsilon$ values. Our numerical results (figures not shown) agree with our analytical work.

\section{Conclusion}

\indent{}\indent{}In the three species RPS system, the enemy of your enemy is \textit{not} your friend. In fact, you are the \textbf{enemy} of your enemy's enemy (just add more ``enemy of the enemy of...'' when you add more species to the system). Said another way, when $N=3$, the enemy of your enemy is your \textbf{victim}. A species' success is, ultimately, its undoing. That is, the growth of one species population allows its enemy to thrive because the growing species drives the enemy's own rival to extinction. The idea that a species' success triggers its own downfall via a negative feedback loop, or ``Survival of the Weakest'' \cite{frean2001rock}, is why we see cyclic dynamics. In contrast, an even number of species enables the community to be divided into even teams (odd parity communities always leave an odd man out). This even division into alliances means that the enemy of your enemy is your \textbf{friend}. A species' population growth only furthers its own success because it helps its friend drive the enemy to extinction. Rather than cycles of competitive exclusion, we see positive feedback loops driving stable alliances. Thus, adding a species to the community effectively ``flips the sign'' of the positive/negative feedback loop.\\
\indent{}This paper was centered on understanding how communities can adapt to and transition between changes in parity. We considered a modified May-Leonard model \cite{may1975nonlinear} with a symmetric inter-species competition coefficient $\alpha$ and introduced a parameter $\varepsilon$ to trace a homotopy between the odd parity rock-paper-scissors (RPS) system and the even parity Hand of Bridge system. We showed that our model exhibits either stable coexistence ($\alpha$ is small) or cycling behavior ($\alpha$ is large) and transitions towards stable alliances of noncompeting species as $\varepsilon$ increases from zero to unity. We also identified transitionary stable steady states (e.g., three-species equilibria and bistability between three- and two-species alliances) absent from both original systems with symmetric competition coefficients and unique to our model.\\
\indent{}We then considered the three different heteroclinic cycles ($H_{34}$, $H_4$, $H_3$) in this system and analyzed their stability. We showed that intermediate values of $\alpha$ ($2<\alpha<3$) allow for large windows in parameter space in which the system maintains bi- and tri- stability between cycles and steady states, while large values of $\alpha$ confine cycling behavior to increasingly extreme values of $\varepsilon$ (i.e., $\varepsilon$ values nearing zero and unity as $\alpha \rightarrow \infty$). That is, when inter-species competition is very high, stable cycling behavior is only possible when species interactions exhibit cyclic competition with nearly perfect symmetry.\\
\indent{}Through numerical simulations, we also identified distinct regimes in parameter space where the type of cycling behavior and the degree to which solution trajectories are attracted to cycles depends on the strength of inter-species competition and the level at which the fourth species participates in the loop of cyclic competition. Our numerical computations showed the existence of stable limit cycles within a window of intermediate $\alpha$ values. Additionally, we found that stable limit cycles and heteroclinic cycles can coexist in our system. A complementary in-depth analysis pertaining to the boundaries of existence and stability of limit cycles is a logical next step for this project. Other studies have shown that an RPS system with mutations and broken symmetries can sustain limit cycles \cite{toupo2015nonlinear,mobilia2010oscillatory,yusuf2026analysis}, and future research should consider how ecological communities transition from odd to even parity when there exists asymmetry in reproductive rates and intra- and inter-species competition.\\
\indent{}Our results demonstrate how seemingly simple ecological communities can display complex behavior that is highly dependent on the number and strength of species interactions.



\section*{Declaration of competing interests}

The authors declare the following financial interests/ personal relationships which may be considered as potential competing interests: TMM reports receiving financial support by the National Science Foundation. The other authors declare that they have no known competing financial interests or personal relationships that could have appeared to influence the work reported in this paper. 

\section*{Acknowledgments}

TMM is supported by the National Science Foundation Graduate Research Fellowship Program under Grant No. DGE-2234667. This research did not receive any specific grant from funding agencies in the public, commercial, or not-for-profit sectors.

\appendix

\section{Analysis of $\mathbf{f_{234}}$}\label{apdx:f234}

The associated Jacobian Matrix is
\begin{align*}
    J_{234} = \begin{bmatrix}
    k_1 & 0 & 0 & 0 \\
    0 & -1+\alpha-\varepsilon\alpha^2  & k_1-1 & 0 \\
    k_2 & 0 & -1+\varepsilon\alpha & -\varepsilon\alpha+\varepsilon^2\alpha^2  \\
    -\alpha & 0 & 0 & -1
\end{bmatrix},
\end{align*}
where 
\begin{align*}
    k_1 = 1-\alpha+\alpha^2-\varepsilon\alpha^3 \text{ and } k_2 = -\alpha+\varepsilon\alpha+\varepsilon\alpha^2-\varepsilon^2\alpha^2.
\end{align*}
The matrix $J_{234}$ has eigenvalues 
\begin{align*}
    &\lambda_1=1-\alpha+\alpha^2-\varepsilon\alpha^3, \lambda_2=-1+\alpha-\varepsilon\alpha^2,\\
    &\lambda_3 = -1+\varepsilon\alpha, \text{ and } \lambda_4=-1.
\end{align*}
When $f_{234}$ is physical, $\lambda_{2},\lambda_{3}<0$.\ We see that $\lambda_{1}<0$ when
\begin{align*}
    \varepsilon_{**}=\frac{1}{\alpha}-\frac{1}{\alpha^2}+\frac{1}{\alpha^3}<\varepsilon.
\end{align*}

\section{Analysis of $\mathbf{f_{123}}$}\label{apdx:f123}
For convenience, we write $f_{123}=(u^*,v^*,w^*,0)$.\ The associated Jacobian Matrix is
\begin{align*}
    J_{123} = \begin{bmatrix}
    -u^* & -\alpha u^* & 0 & 0 \\
    0 & -v^*  & -\alpha v^* & 0 \\
    -(1-\varepsilon)\alpha w^* & 0 & -w^* & -\varepsilon\alpha w^* \\
    0 & 0 & 0 & 1-\alpha u^*
\end{bmatrix}
\end{align*} 
with characteristic equation
\begin{align*}
    P_4(\lambda)=&\ (\lambda-(1-\alpha u^*))\\
    &[(u^*+\lambda)(v^*+\lambda)(w^*+\lambda)+u^*v^*w^*\alpha^3(1-\varepsilon)]=0
\end{align*} and eigenvalues 
\begin{align*}
    &\lambda_1 = -1, \\
    &\lambda_2 = \frac{1}{2}\left[ 1-(u^*+v^*+w^*)+\beta\right],\\
    &\lambda_3 = \frac{1}{2}\left[ 1-(u^*+v^*+w^*)-\beta\right], \\
    &\lambda_4 = 1-\alpha u^*=\frac{1-\alpha+\alpha^2-\varepsilon\alpha^3}{1+\alpha^3(1-\varepsilon)},
\end{align*} 
where 
\begin{align*}
    \beta = \sqrt{(u^*+v^*+w^*-1)^2-4u^*v^*w^*(1+\alpha^3(1-\varepsilon))}.
\end{align*}
We find that $\lambda_4<0$ when 
\begin{align*}
    \varepsilon>\frac{1}{\alpha}-\frac{1}{\alpha^2}+\frac{1}{\alpha^3}=\varepsilon_{**}.
\end{align*}
Additionally, $\lambda_2<0$ and $\lambda_3<0$ when 
\begin{align*}
    \varepsilon>\frac{(\alpha-1)^3-1}{\alpha(\alpha-1)^2}=\varepsilon_{123}.
\end{align*}
Therefore, the state $f_{123}$ is stable when both the following conditions are satisfied:
\begin{align*}
    &\varepsilon_{**}=\frac{1}{\alpha}-\frac{1}{\alpha^2}+\frac{1}{\alpha^3}<\varepsilon \text{ and } \varepsilon_{123}=\frac{(\alpha-1)^3-1}{\alpha(\alpha-1)^2}<\varepsilon.
\end{align*}
We also find that $\lambda_2$ and $\lambda_3$ are purely real when 
\begin{align*}
    \varepsilon>\frac{\left(1-\alpha+\alpha^{2}\right)\left(3\alpha^2-2\alpha+1-2\sqrt{1-\alpha+\alpha^{2}}\right)}{(\alpha-1)(3\alpha^3-\alpha^2+\alpha+1)}.
\end{align*} 
We note that 
\begin{align*}
    &\varepsilon>\frac{\left(1-\alpha+\alpha^{2}\right)\left(3\alpha^2-2\alpha+1-2\sqrt{1-\alpha+\alpha^{2}}\right)}{(\alpha-1)(3\alpha^3-\alpha^2+\alpha+1)} \\
    &\implies \varepsilon>\varepsilon_{123}.
\end{align*}

\section{Analysis of $\mathbf{f_{1234s}}$}\label{apdx:f1234s}
The associated Jacobian matrix is
\begin{align*}
    J_{1234s} = \begin{bmatrix}
    k & \alpha k & 0 & 0 \\
    0 & k & \alpha k & 0 \\
    (1-\varepsilon)\alpha k & 0 & k & \varepsilon \alpha k  \\
    \alpha k & 0 & 0 & k
\end{bmatrix},
\end{align*}
where $k = \frac{-1}{\alpha+1}$. The characteristic equation is 
\begin{align*}
    &P_4(\lambda) = (k-\lambda)^4+(\alpha k)^3(1-\varepsilon)(k-\lambda)-\varepsilon(\alpha k)^4=0. 
\end{align*} 
Let $k-\lambda = \alpha k x$. Then,
\begin{align*}
    &(\alpha k)^4x^4+(\alpha k)^4(1-\varepsilon)x-\varepsilon(\alpha k)^4=0,\\
    &x^4+(1-\varepsilon)x-\varepsilon=0.
\end{align*}
By inspection, we see that $x=-1$ is a root of the above equation.\ Thus,
\begin{align*}
    &\lambda_1 = k(1+\alpha)=-1,\\
    &(x+1)(x^3-x^2+x-\varepsilon)=0.
\end{align*}
Then, the remaining eigenvalues are
\begin{align*}
    \lambda_{j} = \frac{\alpha x_j-1}{\alpha+1},
\end{align*}
where $x_j$ is the $j^{th}$ root of
\begin{align*}
    x^3-x^2+x-\varepsilon=0.
\end{align*}

\section{Analysis of $\mathbf{f_{1234a}}$}\label{apdx:f1234a}
The associated Jacobian matrix is 
\begin{align*}
    &J_{1234a} =\\[1ex] &\begin{bmatrix}
    -U & -\alpha U & 0 & 0 \\
    0 & \frac{U-1}{\alpha}  & U-1 & 0 \\
    (\varepsilon_{**}-1)\left(\frac{\alpha-1+U}{\alpha}\right) & 0 & \frac{1-U-\alpha}{\alpha^2} & \varepsilon_{**}\left(\frac{1-U-\alpha}{\alpha}\right) \\
    \alpha (\alpha U -1) & 0 & 0 &(\alpha U -1) 
\end{bmatrix}
\end{align*}
with eigenvalues
\begin{align*}
    &\lambda_{1} = 0,\\
    &\lambda_{2} = -1,\\
    &\lambda_{3} = \frac{1}{2\alpha^2}(\beta + \sqrt{1-\omega_1+\omega_2}),\\
    &\lambda_{4} = \frac{1}{2\alpha^2}(\beta - \sqrt{1-\omega_1+\omega_2}),
\end{align*}
where
\begin{align*}
        &\beta = 1-2\alpha + U(\alpha^3-\alpha^2+\alpha-1),\\
        &\omega_1 = 2U(1+\alpha-\alpha^2+3\alpha^3-2\alpha^4+2\alpha^5),\\
        &\omega_2 = U^2(1+2\alpha-\alpha^2+4\alpha^3-\alpha^4+2\alpha^5+\alpha^6).
\end{align*}
A basis for the null space of the Jacobian matrix is given by
\begin{align*}
    \text{Null}(J) = \text{span}\left\{ \begin{bmatrix}
    -\frac{1}{\alpha} \\
    \frac{1}{\alpha^2} \\
    -\frac{1}{\alpha^3} \\
    1 
\end{bmatrix} \right\}.
\end{align*}
Thus, we find that the states $f_{1234a}$ lie on the following line:
\begin{align*}
    \begin{bmatrix}
    U \\
    \frac{1-U}{\alpha} \\
    \frac{\alpha-1+U}{\alpha^2} \\
    1-\alpha U 
    \end{bmatrix} + s    \begin{bmatrix}
    -\frac{1}{\alpha} \\
    \frac{1}{\alpha^2} \\
    -\frac{1}{\alpha^3} \\
    1 
    \end{bmatrix}.
\end{align*}
Furthermore, we can derive a conserved quantity $C$ such that $\nabla C \cdot \mathbf{f} = 0$, where $\mathbf{f}$ is a vector describing the system of equations. It follows that
\begin{align*}
    C = \ln u -\alpha\ln v + \alpha^2 \ln w -\varepsilon\alpha^3\ln z=\ln \left( \frac{uw^{\alpha^2}}{v^{\alpha}z^{\alpha^2-\alpha+1}}\right).
\end{align*}
Equivalently, we see that 
\begin{align*}
    K=e^C =  \frac{uw^{\alpha^2}}{v^{\alpha}z^{\alpha^2-\alpha+1}}.
\end{align*}
Each fixed point along the continuous curve has one corresponding zero eigenvalue whose eigenvector points in the direction of the line of steady states. Any trajectory beginning on the line will remain at its initial point, but any perturbation will lead to ``stability with shift'', in which the solution will evolve to a different steady state on the curve such that $K=\text{constant}$ is satisfied.

\section{$H_{34}$ Step 1}\label{apdx:H34}

We begin with the point 
\begin{align*}
    B_1=(\mathcal{O}(1), \delta, \mathcal{O}(1),\mathcal{O}(1)),
\end{align*}
which lies close to the separatrix $f_{34} \rightarrow f_1$. We linearize (\ref{eqs:3to4_system}) about $f_1$ by considering the change of variables
\begin{align*}
    u \sim 1+\tilde{u}, v \sim \tilde{v}, w\sim \tilde{w}, \text{ and } z\sim\tilde{z}.
\end{align*}
We then obtain
\begin{align*}
    &\tilde{u}'=-\tilde{u}-\alpha\tilde{v},\\
    &\tilde{v}'=\tilde{v},\\
    &\tilde{w}'=[1-(1-\varepsilon)\alpha]\tilde{w},\\
    &\tilde{z}'=(1-\alpha)\tilde{z}.
\end{align*}
We solve the linearized system with $B_1$ as the initial conditions, i.e.,
\begin{align*}
    \tilde{u}(0)=u_1-1, \tilde{v}(0)=\delta, \tilde{w}(0)=w_1, \text{ and } \tilde{z}(0)=z_1,
\end{align*}
and get
\begin{align*}
    &\tilde{u}(t)=(u_1-1+\frac{1}{2}\alpha\delta)e^{-t}-\frac{1}{2}\alpha\delta e^t,\\
    &\tilde{v}(t)=\delta e^t,\\
    &\tilde{w}(t)=w_1 e^{-[(1-\varepsilon)\alpha-1]t},\\
    &\tilde{z}(t)=z_1 e^{-(\alpha-1)t}.
\end{align*}
Let $t_1$ be the time at which the trajectory reaches the point $B_2$. We can determine $t_1$ by solving for the time at which $\tilde{v}(t)$ grows from the small coordinate $\delta$ at $B_1$ to the $\mathcal{O}(1)$ quantity $v_2$ at $B_2$. Thus,
\begin{align*}
    \tilde{v}(t_1)=v_2 &\implies \delta e^{t_1}=v_2 \\&\implies e^{t_1}=v_2\delta^{-1} \\
    &\implies t_1 = \ln{(v_2)} - \ln{(\delta)}.
\end{align*}
We then substitute $t=t_1$ into the solutions found above to find the remaining coordinates of the point $B_2$:
\begin{align*}
    &u_2 = u(t_1)=1+\tilde{u}(t_1)=1+(u_1-1+\frac{1}{2}\alpha\delta)\frac{\delta}{v_2}-\frac{1}{2}\alpha v_2,\\
    &w_2 = \tilde{w}(t_1)=w_1 (\frac{\delta}{v_2})^{(1-\varepsilon)\alpha-1},\\
    &z_2 = \tilde{z}(t_1)=z_1 (\frac{\delta}{v_2})^{\alpha-1}.
\end{align*}
We are only interested in tracking the small coordinates of each point around the cycle. The point $B_2$ lies near the separatrix $f_1 \rightarrow f_2$, so the coordinates $w_2$ and $z_2$ must be small. Thus, we impose the conditions
\begin{align*}
    \varepsilon<1-\frac{1}{\alpha} \text{ and } \alpha>1
\end{align*}
and conclude
\begin{align*}
    &B_1=(\mathcal{O}(1), \delta, \mathcal{O}(1),\mathcal{O}(1)) \rightarrow \\
    &B_2 = (\mathcal{O}(1),\mathcal{O}(1),\mathcal{O}(1) \delta^{(1-\varepsilon)\alpha-1},\mathcal{O}(1)\delta^{\alpha-1}).
\end{align*}

\section{$H_4$ is unstable when $\lambda_{1,2}$ are complex conjugates}\label{apdx:H4}

We remind the reader that $H_4$ is stable when the components of $\vec{V}_n$ tend to negative infinity as $n \rightarrow \infty$, where
\renewcommand{\arraystretch}{1.0}{
\begin{align*}
    \vec{V}_n = \mathbf{A}^n\vec{V}_{0} \text{ and } \mathbf{A} = \begin{bmatrix}
        a & b\\
        c & d
    \end{bmatrix}.
\end{align*}}
We consider the case in which the matrix $\mathbf{A}$ has complex conjugate eigenvalues $\lambda_1$ and $\lambda_2$. We replace $\lambda_1$ with 
$\lambda=r e^{i\theta}$ and $\lambda_2$ with $\bar{\lambda}=r e^{-i\theta}$ 
(we can always assume that $0<\theta<\pi$). After diagonalizing $\mathbf{A}$ and determining $\mathbf{A}^n$, we get
\begin{align*}
    \vec{V}_n &= \frac{1}{c(\lambda-\bar{\lambda})}
{\begin{bmatrix} 
v_1\\[1ex]
v_2
\end{bmatrix}} =
\frac{r^{n-1}}{\sin\theta}
{\begin{bmatrix} 
\omega_1 \\[1ex]
\omega_2
\end{bmatrix}},
\end{align*}
where
\begin{align*}
    v_1= &\ c[-(d-\lambda)\lambda^n + (d-\bar{\lambda})\bar{\lambda}^n]\ln (\delta) -\\
&(d-\lambda)(d-\bar{\lambda})(\lambda^n-\bar{\lambda}^n) \ln (\sigma),\\
v_2 = &\ c^2 (\lambda^n-\bar{\lambda}^n) \ln(\delta) +
c [ \lambda^n (d-\bar{\lambda})- \bar{\lambda}^n (d-\lambda)]\ln(\sigma), \\
\omega_1 = &\ [-d \sin(n\theta) + r \sin(n+1)\theta] \ln (\delta) -\\
&\frac{1}{c} \sin(n\theta) (d^2-2dr\cos(n\theta)+r^2) \ln (\sigma),\\
\omega_2 =&\ c \sin(n\theta)  \ln(\delta) +
 [ d\sin(n\theta) -r \sin(n-1)\theta ]\ln(\sigma).
\end{align*}
We now prove that $\vec{V}_n$ does not go to $-\infty$ as $n\to\infty$ and, therefore,
the HH4-cycle is unstable in the case of complex eigenvalues.

Let us consider separately two cases: when $\theta$ and $\pi$
are commensurable and when they are not. If they are not commensurable then,
according to the Kronecker density theorem, there are infinitely many values of $n$
such that $\sin(n\theta)$ is as close to zero as we want while $\sin((n+1)\theta)<0$.
Then, neglecting $\sin(n\theta)$, we see that the first row of $\vec{V}_n$ reduces to 
\begin{align*}
    \frac{r^n}{\sin(\theta)}  \sin((n+1)\theta) \ln (\delta)>0,
\end{align*}
i.e., has a wrong sign for stability.

To analyze the case of commensurable $\theta$ and $\pi$, we will need the following
summation formulas:
\begin{align*}
    \sum_{n=M}^{N} &\sin(n\theta) = \\
    &\frac{1}{\sin(\frac{1}{2}\theta)} 
\sin \left (\frac{1}{2}(N-M+1)\theta \right )\sin \left (\frac{1}{2}(N+M)\theta \right )
\end{align*}
and
\begin{align*}
    \sum_{n=M}^{N} &\sin((n-1)\theta) = \\
    &\frac{1}{\sin(\frac{1}{2}\theta)} 
\sin \left (\frac{1}{2}(N-M+1)\theta \right )\sin \left (\frac{1}{2}(N+M-2)\theta \right ).
\end{align*}
 If $H_4$ were stable, we would have
\begin{align}\label{eqs:H_4_stable_complex}
    c \sin(n\theta)  \ln(\delta) + [ d\sin(n\theta) -r \sin(n-1)\theta ]\ln(\sigma) <0,
\end{align}
for all sufficiently large $n$, and the sum of such inequalities for these $n$
would also be negative. However, since $\theta$ and $\pi$
are commensurable, we have $\theta=\frac{p}{q}\pi$, where $p$ and $q$ are integers.
Then if 
\begin{align*}
     \frac{1}{2}(N-M+1)
\end{align*}
is any multiple of $q$, the sums of the left-hand side of (\ref{eqs:H_4_stable_complex}) for sufficiently large $n$ are equal to zero. Therefore, condition (\ref{eqs:H_4_stable_complex}) is not satisfied and $H_4$ is not stable.

\printcredits

\bibliographystyle{elsarticle-num.bst}

\bibliography{cas-refs}



\end{document}